# FREE PRODUCTS OF SEMISIMPLE TENSOR CATEGORIES

SHIGERU YAMAGAMI

Department of Mathematical Sciences
Ibaraki University
Mito, 310-8512, JAPAN
e-mail: yamagami@mito.ipc.ibaraki.ac.jp

Abstract. Free products of semisimple tensor categories are constructed with the help of polygonal presentation.

## Introduction

In this paper, we shall construct free products of semisimple tensor categories.

A similar notion for subfactors is introduced by D. Bisch and V. Jones [1, 2] and has been a source of their planar algebras (see [9] for more information on planar algebras). The notion turns out to be of combinatorial nature and the construction of free products are worked out for paragroups (a good combinatorial invaraint of subfactors) by S. Gnerre [5]. Although the way of definition is fairly obvious, it involves a lot of case-studies and a complete description of construction is unexpectedly messy without methods. Moreover, these approaches use S. Popa's realization theorem of subfactors ([15]) in a crucial way and seem to be difficult to be extented to tensor categories.

Compared with the precedents, our approach is purely combinatorial and can be generalized to more algebraic situations.

Let us begin with a brief discussion on combinatorial aspects of tensor categories; more precisely we shall work with polygonal presentations of semisimple tensor categories explored in [23]: tensor categories are constructed in terms of associativity transformations among tensor products of triangular vector spaces indexed by triplets in the fusion rule set, where the associativity transformations are required to satisfy the pentagonal relations.

The free product of fusion algebras (or the fusion rule set) is easy to define and is again a fusion algebra as seen by Hiai and Izumi in [8]. For triangular vector spaces as well as associativity transformations, it is also not difficult to introduce their free products by appealing to inductive definitions.

The proof of pentagonal relations is, however, unduly complicated with a large number of case-classifications if one tries to check them in a naive way. To remedy such difficulties, we apply a method obtained during studies of polygonal coherence theorem to polygonal vector spaces of free products, which enables us to reduce the whole description of relevant cases to a moderate amount.

When tensor categories are restricted to more specific classes such as rigid or C*-tensor categories, our construction immediately provides the corresponding versions of free products too. In particular, free products of rigid C*-tensor categories remain in the same class. Since tensor categories in this class are exactly the ones realized by von Neumann algebra (more precisely factor) bimodules of finite Jones index ([22]), we can even recover the results for subfactors or paragroups as well.

As another application of the present construction, we shall identify planar algebras of Bisch and Jones as the ascending algebra associated to the $m$-times free product of the fundamental generator in the Tannaka dual of the quantum group $SL_q(2,\mathbb{C})$, which particularly gives the semisimplicity criterion to the planar algebras: the Fuss-Catalan algebras $FC_n(a_1,\dots,a_m)$ with the evaluation parameters $a_1,\dots,a_m$ are semisimple for all $n \geq 1$ if and only if none of $a_1,\dots,a_m$ belongs to the set

$$\{2\cos\pi r | r \in \mathbb{Q}\setminus\mathbb{Z}\}.$$

Notation and Terminology

By a semisimple tensor category, we shall mean a tensor category based on the complex number field $\mathbb{C}$ such that the unit object is simple and any object is isomorphic to a direct sum of finitely many simple objects.

Given a semisimple tensor category $\mathcal{C}$, the set of isomorphism classes of simple objects is denoted by $\mathrm{Spec}(\mathcal{C})$ and is referred to as the spectrum set. The spetrum set is often identified with its representative set $S$ consisting of simple objects in $\mathcal{C}$ including the unit object $I$ (or 1).

Given such a representative set $S$ of $\mathrm{Spec}(\mathcal{C})$, we write

$$\mathrm{Hom}(x\otimes y, z) = \begin{bmatrix} x\,y \\ z \end{bmatrix} = \underset{z}{x \bigtriangleup y}$$

for $x$, $y$, $z \in S$ and call it a triangular vector space.

Given $x_j \in S$ with $0 \leq j \leq 3$, the vector space $\mathrm{Hom}(x_1 \otimes x_2 \otimes x_3, x_0)$ is decomposed into a direct sum of triangular vector spaces in two ways: (i) first decompose $x_1 \otimes x_2$ in terms of $x_{12} \in S$ and then take out the $x_0$-component in $x_{12} \otimes x_3$ or (ii) first decompose $x_2 \otimes x_3$ in terms of $x_{23} \in S$ and then take out the $x_0$-component in $x_1 \otimes x_{23}$.

These are denoted as

$$x_1 \overset{x_2}{\underset{x_0}{\boxslash}} x_3 = \bigoplus_{x_{12} \in S} \begin{bmatrix} x_1 \ x_2 \\ x_{12} \end{bmatrix} \otimes \begin{bmatrix} x_{12} \ x_3 \\ x_0 \end{bmatrix},$$

$$x_1 \overset{x_2}{\underset{x_0}{\boxbslash}} x_3 = \bigoplus_{x_{23} \in S} \begin{bmatrix} x_2 \ x_3 \\ x_{23} \end{bmatrix} \otimes \begin{bmatrix} x_1 \ x_{23} \\ x_0 \end{bmatrix}.$$

The associativity transformation $a_{x_0,x_1,x_2,a_3}$ is then defined by connecting these decompositions:

$$a_{x_0,x_1,x_2,a_3} : x_1 \overset{x_2}{\underset{x_0}{\boxslash}} x_3 \to x_1 \overset{x_2}{\underset{x_0}{\boxbslash}} x_3.$$

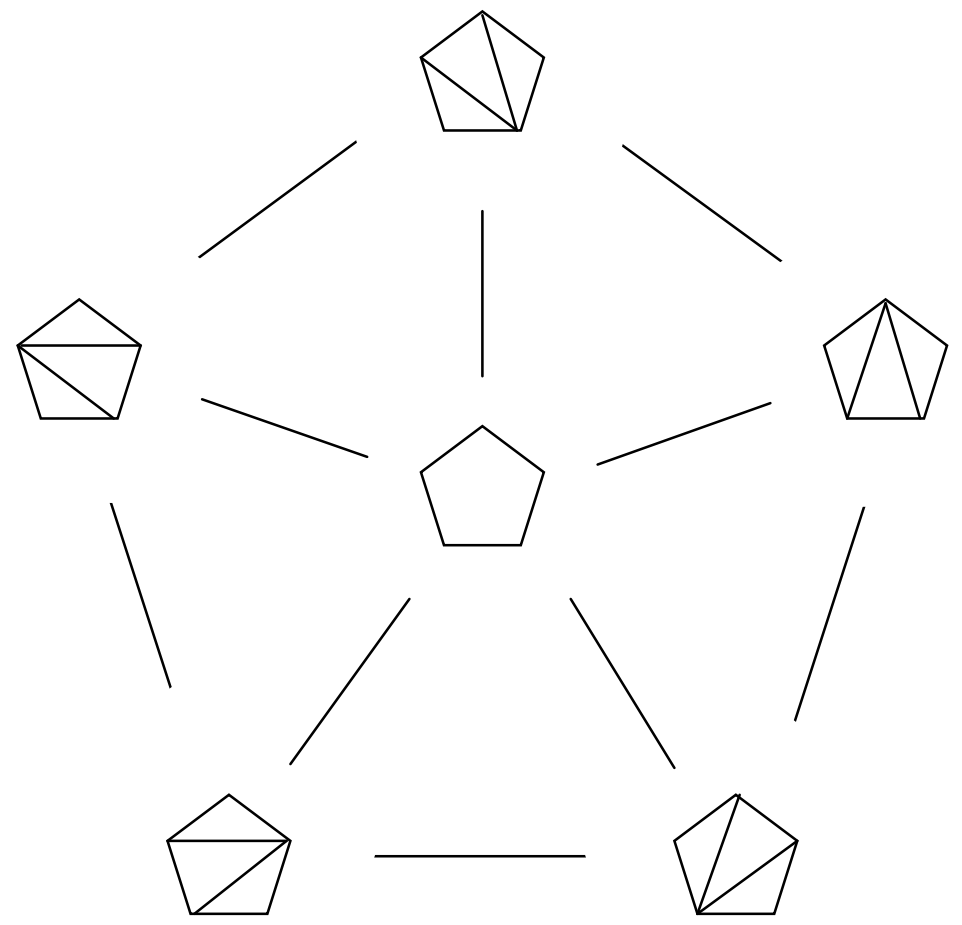

FIGURE 1.

The left and right unit constraints give rise to special vectors

$$l_s \in \begin{bmatrix} 1\, s \\ s \end{bmatrix}, \qquad r_s \in \begin{bmatrix} s\, 1 \\ s \end{bmatrix},$$

which satisfy the unit constraint condition

$$a(r_x \otimes \zeta) = l_y \otimes \zeta \qquad \text{for } \zeta \in \begin{bmatrix} x\, y \\ z \end{bmatrix}$$

with $a = a_{z,x,1,y}$.

For the pentagonal vector space $\mathrm{Hom}(x_1 \otimes x_2 \otimes x_3 \otimes x_4, x_0)$, we have five ways of decompositions, which are connected by (amplified) associativity transformations as Fig. 1. The outcircuit diagram is then commutative, which is referred to as the pentagonal relation for the associativity transformations.

Conversely, starting with (triangular) vector spaces indexed by triplets in a set $S$ with a distinguished element 1, a family of non-zero vectors $\{l_s, r_s\}_{s\in S}$ and a family of isomorphisms $\{a_{x_0,x_1,x_2,x_3}\}$ fulfilling the unit constraint condition and the pentagonal relation, we can reconstruct the tensor category (see [16, 23]).

## 1. Associativity Transformations

Let $\mathcal{C}$ and $\mathcal{D}$ be semisimple tensor categories with spectrum (fusion rule) sets $S$ and $T$ respectively. Let $S * T$ be the free product of $S$ and $T$ on which fusion rule is defined in a free fashion. Note that

$$S * T = \{1\} \sqcup S^\circ \sqcup T^\circ \sqcup (S^\circ \times T^\circ) \sqcup (T^\circ \times S^\circ) \sqcup (S^\circ \times T^\circ \times S^\circ) \sqcup (T^\circ \times S^\circ \times T^\circ) \sqcup \dots$$

with

$$S^\circ = S \setminus \{1\}, \qquad T^\circ = T \setminus \{1\}.$$

In what follows, we use the following notation: for $x$, $y \in S * T$, we write $x \| y$ or $x|y$ according to $xy \in S * T$ with $|xy| = |x| + |y|$ or not ($|x|$ denotes the word length of $x$ with the convention $|1| = 0$).

Given a triplet $(x, y, z)$ in $S * T$, we introduce the associated triangular vector space in the following way: For $x\|y$, we set

$$\begin{bmatrix} x\, y \\ z \end{bmatrix} = \mathbb{C}\delta_{xy,z}.$$

(Note here that this includes the case that $x = 1$ or $y = 1$.) If $x|y$, we define the triangular vector space by induction on the length $|x| + |y|$ ($|x| + |y| \geq 2$ by assumption). For the lowest case $|x| = 1 = |y|$, we set

$$\begin{bmatrix} x\, y \\ z \end{bmatrix} = \begin{cases} \mathrm{Hom}(x \otimes y, z) & \text{if } x,\ y \text{ and } z \text{ have the same parity,} \\ \{0\} & \text{otherwise.} \end{cases}$$

For $|x|+|y|\geq 3$ with $x=x'a$ and $y=by'$ ($|a|=|b|=1$), we require

$$\begin{bmatrix} x\,y \\ z \end{bmatrix} = \begin{cases} \begin{bmatrix} a\,b \\ c \end{bmatrix} & \text{if } z=x'cy' \text{ with } |c|=1, \\ \begin{bmatrix} a\,b \\ 1 \end{bmatrix} \otimes \begin{bmatrix} x'\ y' \\ z \end{bmatrix} & \text{otherwise.} \end{cases}$$

We can easily see the following.

**Lemma 1.1** (Squeeze). *The triangular vector space* $\begin{bmatrix} x\,y \\ z \end{bmatrix}$ *is non-trivial only if* $z=x''wy''$ *with* $x=x''x'$, $y=y'y''$, $|x'|=|y'|$ *and* $|w|=1$ *or* $z=x\ominus y$, *i.e., only if* $z$ *is a 'squeeze' of* $xy$.

We next introduce associativity transformations

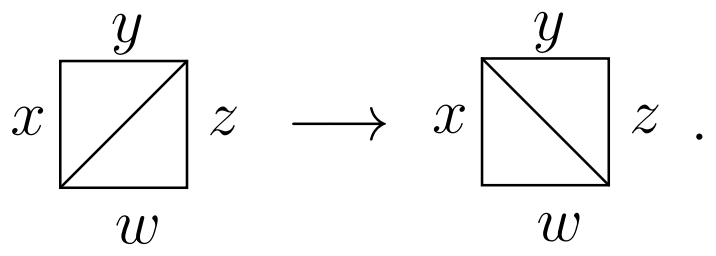

For $x\|y$ and $y\|z$, we have

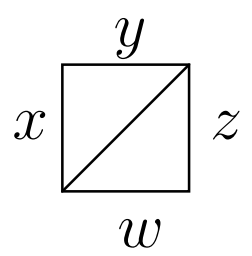

$= \mathbb{C}\,\delta_{xyz,w} =$
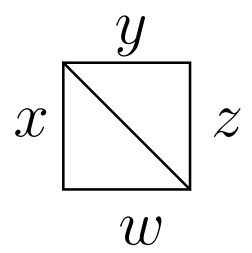

If $x\|y$ and $y|z$, the transformation is defined inductively on the length $|y|$: For $|y|=1$, letting $z=cz'$, we have

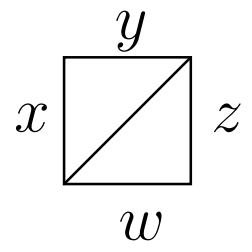

$$= \bigoplus_{u\neq 1} \delta_{xuz',w} \begin{bmatrix} y\,c \\ u \end{bmatrix} \oplus \begin{bmatrix} y\,c \\ 1 \end{bmatrix} \otimes \begin{bmatrix} x\,z' \\ w \end{bmatrix} =$$
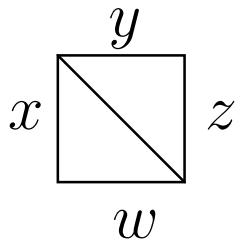

For $|y| \geq 2$, letting $y = y'b$, we have

$$x \underset{w}{\overset{y}{\boxed{\diagup}}} z = \bigoplus_{u \neq 1} \delta_{xy'uz',w} \begin{bmatrix} b\,c \\ u \end{bmatrix} \oplus \begin{bmatrix} b\,c \\ 1 \end{bmatrix} \otimes x \underset{w}{\overset{y'}{\boxed{\diagup}}} z'$$

$$x \underset{w}{\overset{y}{\boxed{\diagdown}}} z = \bigoplus_{u \neq 1} \delta_{xy'uz',w} \begin{bmatrix} b\,c \\ u \end{bmatrix} \oplus \begin{bmatrix} b\,c \\ 1 \end{bmatrix} \otimes x \underset{w}{\overset{y'}{\boxed{\diagdown}}} z'$$

and the associativity transformation is defined by the identity plus the ampliation of the isomorphism

$$x \underset{w}{\overset{y'}{\boxed{\diagup}}} z' \longrightarrow x \underset{w}{\overset{y'}{\boxed{\diagdown}}} z',$$

which is well-defined by induction hypothesis.

Let $y = y_k \dots y_1$, $z = z_1 \dots z_l$ with $|y| = k$, $|z| = l$ and write $y(j) = y_k \dots y_{j+1}$, $z(j) = z_{j+1} \dots z_l$ for $j = 1, 2, \dots$.

If we apply the above inductive procedure repeatedly, we obtain the explicit formula:

$$\begin{aligned} x \underset{w}{\overset{y}{\boxed{\diagup}^{\bullet}}} z = \begin{bmatrix} xy\ z \\ w \end{bmatrix} &= \sum_{u_1 \neq 1} \delta_{xy(1)u_1z(1),w} \begin{bmatrix} y_1\ z_1 \\ u_1 \end{bmatrix} + \sum_{u_2 \neq 1} \delta_{xy(2)u_2z(2),w} \begin{bmatrix} y_1\ z_1 \\ 1 \end{bmatrix} \otimes \begin{bmatrix} y_2\ z_2 \\ u_2 \end{bmatrix} \\ &+ \sum_{u_3 \neq 1} \delta_{xy(2)u_3z(3),w} \begin{bmatrix} y_1\ z_1 \\ 1 \end{bmatrix} \otimes \begin{bmatrix} y_2\ z_2 \\ 1 \end{bmatrix} \otimes \begin{bmatrix} y_3\ z_3 \\ u_3 \end{bmatrix} + \dots \\ &+ [y \cap z] \otimes x \underset{w}{\overset{\circ}{\triangle}} {}^{y \triangle z}, \end{aligned}$$

where the vector space $[y \cap z]$ is defined by

$$[y \cap z] = \begin{bmatrix} y_1\ z_1 \\ 1 \end{bmatrix} \otimes \begin{bmatrix} y_2\ z_2 \\ 1 \end{bmatrix} \otimes \dots$$

and

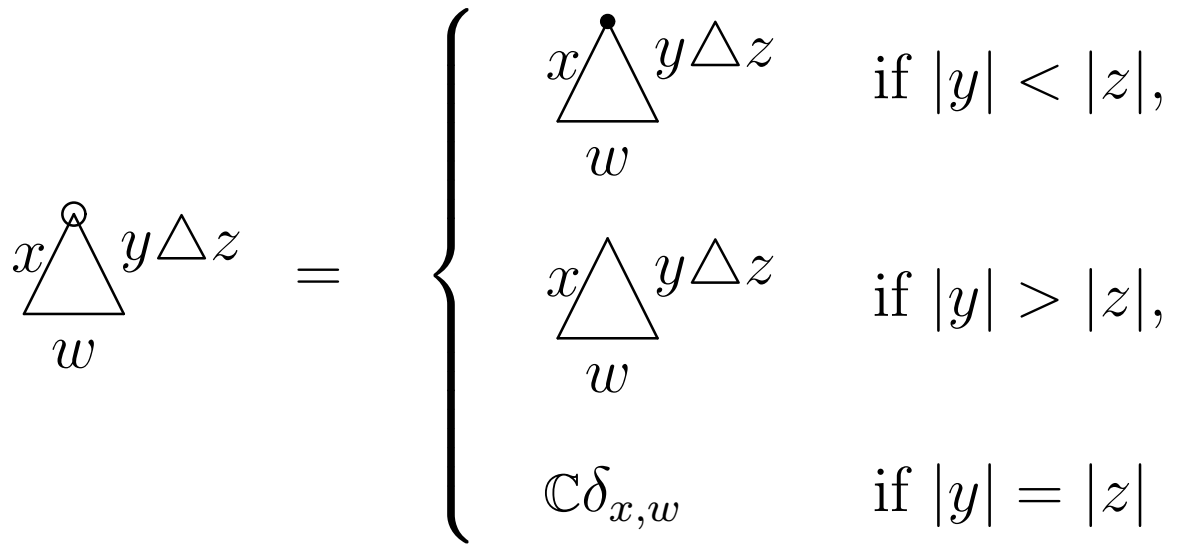

with $y\triangle z \in S * T$ defined by

$$y\triangle z \quad = \begin{cases} \text{the sequence of first } |z|-|y| \text{ elements in } z & \text{if } |y|<|z|, \\ \text{the sequence of first } |y|-|z| \text{ elements in } y & \text{if } |y|>|z|, \\ 1 & \text{if } |y|=|z|. \end{cases}$$

We also use the notation

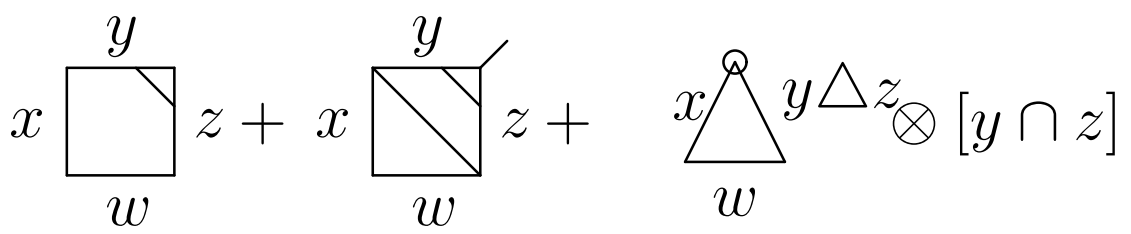

to express the above formula for 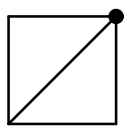 :

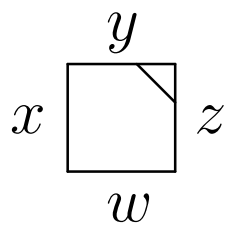 $$= \sum_{u_1\neq 1} \delta_{xy(1)u_1z(1),w} \begin{bmatrix} y_1\ z_1 \\ u_1 \end{bmatrix}$$

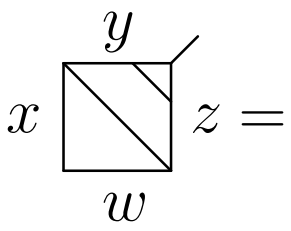 $$= \sum_{u_2\neq 1} \delta_{xy(2)u_2z(2),w} \begin{bmatrix} y_1\ z_1 \\ 1 \end{bmatrix} \otimes \begin{bmatrix} y_2\ z_2 \\ u_2 \end{bmatrix} + \sum_{u_3\neq 1} \delta_{xy(2)u_3z(3),w} \begin{bmatrix} y_1\ z_1 \\ 1 \end{bmatrix} \otimes \begin{bmatrix} y_2\ z_2 \\ 1 \end{bmatrix} \otimes \begin{bmatrix} y_3\ z_3 \\ u_3 \end{bmatrix} + \cdot$$

Likewise we have

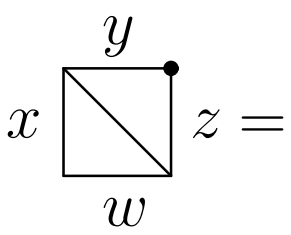 $$= \sum_{u} \begin{bmatrix} y\ z \\ u \end{bmatrix} \otimes \begin{bmatrix} x\ u \\ w \end{bmatrix}$$

$$= \sum_{u} \Big( \sum_{u_1\neq 1} \delta_{y(1)u_1z(1),u} \begin{bmatrix} y_1\ z_1 \\ u_1 \end{bmatrix} \otimes \begin{bmatrix} x\ u \\ w \end{bmatrix}$$

$$+ \sum_{u_2\neq 1} \delta_{y(2)u_2z(2),u} \begin{bmatrix} y_1\ z_1 \\ 1 \end{bmatrix} \otimes \begin{bmatrix} y_2\ z_2 \\ u_2 \end{bmatrix} \otimes \begin{bmatrix} x\ u \\ w \end{bmatrix} + \cdots + \delta_{y\triangle z,u}\, [y\cap z] \otimes \begin{bmatrix} x\ u \\ w \end{bmatrix} \Big)$$

$=$ 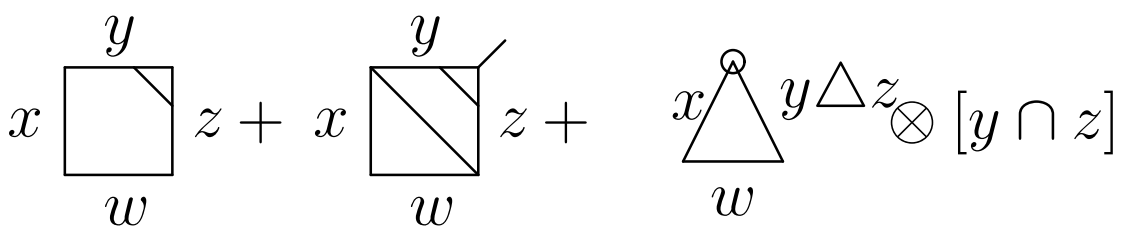

and the asscoaitivity transformation 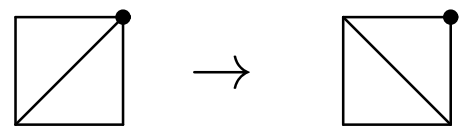 turns out to be the identity.

Similarly we define the transformation if $x|y$ and $y\|z$.

**Lemma 1.2** (Associativity Transformation Formula). *Let $x$, $y$, $z$ and $w \in S * T$ with $x\|y$ or $y\|z$. Then the associativity transformation*

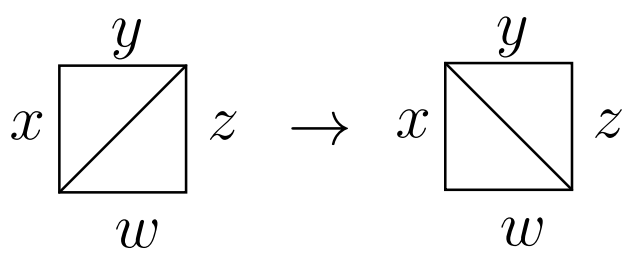

*is the identity when*
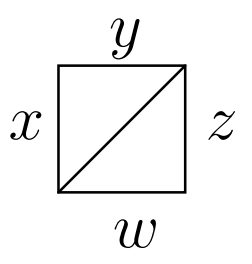

*and*
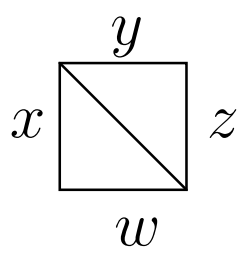

*are expressed by triangulated vector spaces.*

Before going into the remaining case of $x|y$ and $y|z$, we here make a detour to establish a coherence result on associativity transformations defined so far.

Let $x$, $y \in (S * T)^\circ$ be such that $x|y$ and consider a polygon $P$ of $|x|+|y|+1$ edges, labeled by the sequence $x, y$ with an element $z \in S*T$ assigned to the bottom edge. Given a triangulation $\mathcal{T}$ of $P$ and a polygon $Q$ consisting of edges in $\mathcal{T}$, the subspace of the polygonal vector space $[\mathcal{T}]$ specified by a given peripheral labeling $\mathcal{L}$ of $Q$ is non-trivial only when the reduced $Q$ admits at most one vertex of interaction for $\mathcal{L}$. Here the reduced $Q$ means the polygon with unit-assigned edges shrinked to one point.

In fact, letting $\bullet$ be the vertex shared by $x$ and $y$, if $Q$ contains $\bullet$, then the possible labeling $\mathcal{L}$ for non-trivial subspace is just a division of the sequence $x, y$ and hence it has the unique vertex of interaction, i.e., $\bullet$. If $Q$ does not contain $\bullet$, then there is the unique edge $e$ which separates $\bullet$ with the bottom edge of $P$. It is immediate to see that the labeling of other edges is specified by consecutive intervals in $x$ or $y$ for non-trivial subspaces. On the other hand, the labeling of the edge $e$ is forced to be of the form described in Lemma 1.1. There is no interaction if the label $\mathcal{L}(e)$ at $e$ is half squeezed, whereas the interaction occurs for the reduced $Q$ at one of end points of $e$ if $\mathcal{L}(e) = x' \ominus y'$, which may happen to be the unit 1, (Fig. 2).

In particular, the relevant associativity transformations connecting triangulated vector spaces for different choices of $\mathcal{T}$ are those already defined and we can talk about the coherence of them at this stage.

**Lemma 1.3** (One-Vertex Coherence). *Let $x$, $y \in (S * T)^\circ$ and consider a polygon $P$ of $|x| + |y| + 1$-edges labeled by the sequence $x, y$ with $z \in S * T$ assigned to the bottom edge.*

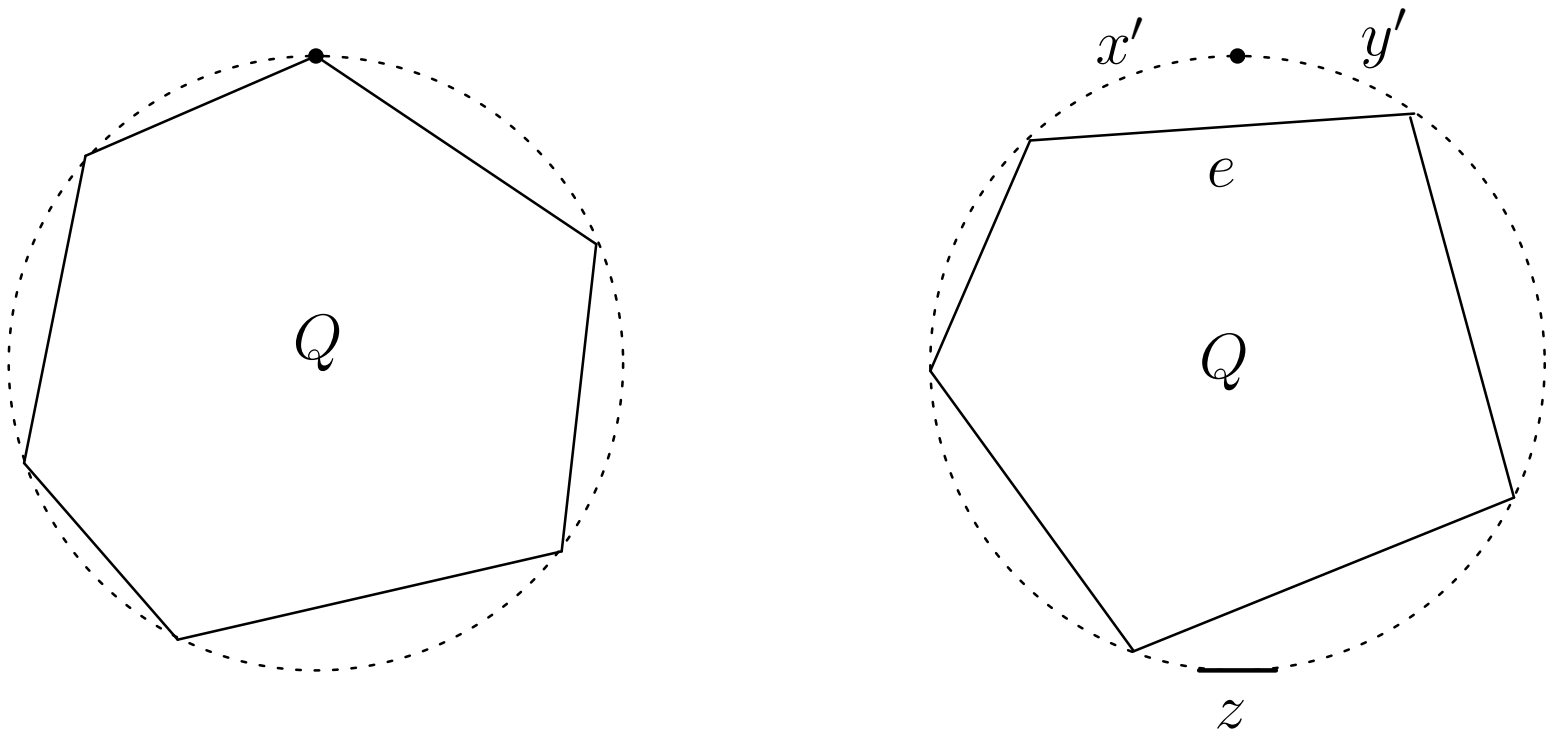

FIGURE 2.

*Then the coherence holds among triangulated vector spaces for possible choices of triangulations of $P$.*

*Proof.* By the proof of coherence theorem in appendix, together with the observation discussed above, the problem is reduced to checking the compatibility of triangulated vector spaces of labeled pentagons which have at most one vertex of interaction.

Consider the pentagon of type $\hom(s \otimes t \otimes t \otimes u \otimes v, w)$ with $s$, $t$, $u$, $v$, $w \in S * T$. If one of $s$, $t$, $u$ and $v$ is equal to the unit element 1, then the pentagonal relation is reduced to the trivial identity of associativity transformations such as

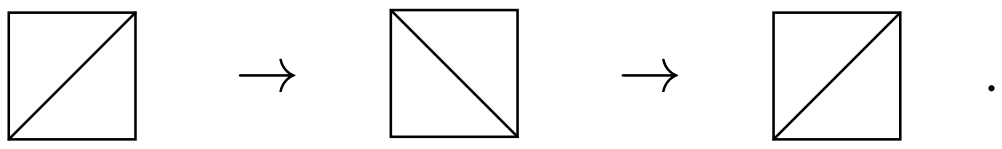

So we may assume that $s$, $t$, $u$ and $v$ belong to $(S * T)^\circ$. If they have no interaction point, the pentagonal relation is trivially satisfied. Otherwise, there are three cases according to the position of the interacting point $\bullet$. Suppose, for example, that the point $\bullet$ is shared by $t$ and $u$. Write $t = t'a$ and $u = bu'$. Then we have

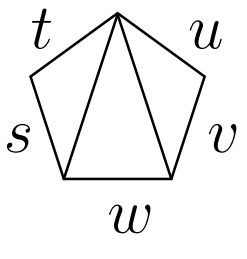

$$= \quad$$

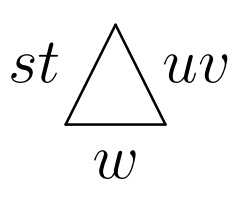

$$= \bigoplus_{c \neq 1} \delta_{st'cu'v,w} \begin{bmatrix} a\,b \\ c \end{bmatrix} + \begin{bmatrix} a\,b \\ 1 \end{bmatrix} \otimes$$

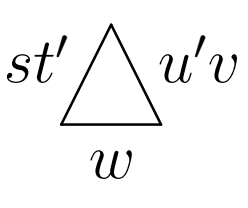

$$= \bigoplus_{c \neq 1} \delta_{st'cu'v,w} \begin{bmatrix} a\,b \\ c \end{bmatrix} + \begin{bmatrix} a\,b \\ 1 \end{bmatrix} \otimes$$

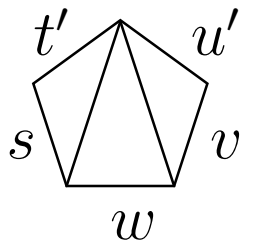

,

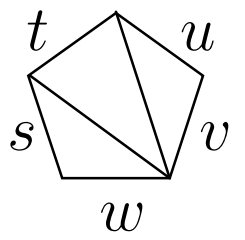 $=$ 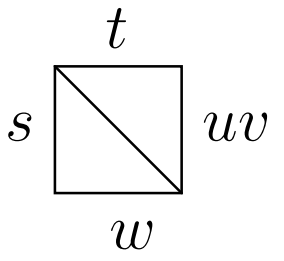 $= \bigoplus_{c\neq 1} \delta_{st'cu'v,w} \begin{bmatrix} a\,b \\ c \end{bmatrix} + \begin{bmatrix} a\,b \\ 1 \end{bmatrix} \otimes$ 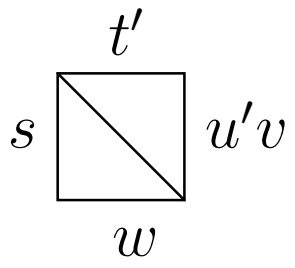

$= \bigoplus_{c\neq 1} \delta_{st'cu'v,w} \begin{bmatrix} a\,b \\ c \end{bmatrix} + \begin{bmatrix} a\,b \\ 1 \end{bmatrix} \otimes$ 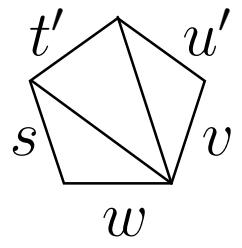

and similarly for others. Thus the coherence is reduced to that for the labeling of lower level $\mathrm{Hom}(s \otimes t' \otimes u' \otimes v, w)$ and we can apply the induction hypothesis. □

We now return to the definition of associativity transformations for $x|y$ and $y|z$. Again we shall appeal to an induction on the length $|y|$.

For $|y| = 1$ with $x = x'a$ and $z = cz'$, we have

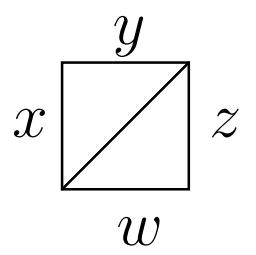 $= \sum_{u\neq 1} \begin{bmatrix} a\,y \\ u \end{bmatrix} \otimes \begin{bmatrix} x'u\ z \\ w \end{bmatrix} + \begin{bmatrix} a\,y \\ 1 \end{bmatrix} \otimes \begin{bmatrix} x'\ z \\ w \end{bmatrix}$

$$= \sum_{\substack{u\neq 1\\ v\neq 1}} \delta_{x'vz',w} \begin{bmatrix} a\,y \\ u \end{bmatrix} \otimes \begin{bmatrix} u\,c \\ v \end{bmatrix} + \begin{bmatrix} a\,y \\ 1 \end{bmatrix} \otimes \begin{bmatrix} x'\ z \\ w \end{bmatrix} + \sum_{u\neq 1} \begin{bmatrix} a\,y \\ u \end{bmatrix} \otimes \begin{bmatrix} u\,c \\ 1 \end{bmatrix} \otimes \begin{bmatrix} x'\ z' \\ w \end{bmatrix}$$

$$= \sum_{v\neq 1} \sum_{u} \begin{bmatrix} a\,y \\ u \end{bmatrix} \otimes \begin{bmatrix} u\,c \\ v \end{bmatrix} \otimes \begin{bmatrix} x'\ vz' \\ w \end{bmatrix} + \sum_{u\neq 1} \begin{bmatrix} a\,y \\ u \end{bmatrix} \otimes \begin{bmatrix} u\,c \\ 1 \end{bmatrix} \otimes \begin{bmatrix} x'\ z' \\ w \end{bmatrix}$$

$= \sum_{v\neq 1}$ 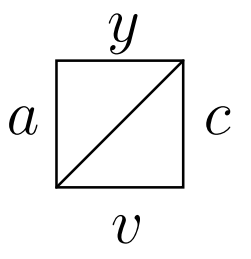 $\otimes \begin{bmatrix} x'\ vz' \\ w \end{bmatrix} + \sum_{u\neq 1} \begin{bmatrix} a\,y \\ u \end{bmatrix} \otimes \begin{bmatrix} u\,c \\ 1 \end{bmatrix} \otimes \begin{bmatrix} x'\ z' \\ w \end{bmatrix}$

$= \sum_{v\neq 1} \delta_{x'vz',w}$ 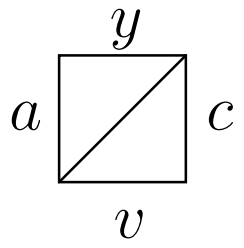 $+$ 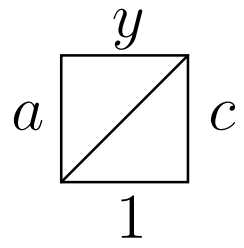 $\otimes \begin{bmatrix} x'\ z' \\ w \end{bmatrix}.$

In the last line, we have used

$$\begin{bmatrix} x'\ vz' \\ w \end{bmatrix} = \mathbb{C}\, \delta_{x'vz',w} \qquad \text{and} \qquad \begin{bmatrix} u\,c \\ 1 \end{bmatrix} = 0 \quad \text{for} \quad u = 1.$$

By symmetry, we have

$$x\,\underset{w}{\overset{y}{\boxbslash}}\,z = \sum_{v\neq 1} \delta_{x'vz',w}\ a\,\underset{v}{\overset{y}{\boxbslash}}\,c \quad + \quad a\,\underset{1}{\overset{y}{\boxbslash}}\,c \otimes \begin{bmatrix} x'\ z' \\ w \end{bmatrix}.$$

Now the associativity transformation for the square $x\,\underset{w}{\overset{y}{\square}}\,z$ with $|y| = 1$ is defined by the help of the isomorphism

$$a\,\underset{v}{\overset{y}{\boxslash}}\,c \longrightarrow c\,\underset{v}{\overset{y}{\boxbslash}}\,a \quad \text{for } v \neq 1 \text{ or } v = 1,$$

which are associativity transformations in $\mathcal{C}$ or $\mathcal{D}$.

For $|y| \geq 2$, writing $x = x'a$, $y = by' = by''b'$ and $z = cz'$, we have

$$x\,\underset{w}{\overset{y}{\boxslash}}\,z = \bigoplus_{u\neq 1} \begin{bmatrix} a\,b \\ u \end{bmatrix} \otimes \begin{bmatrix} x'uy'\ z \\ w \end{bmatrix} \quad + \quad \begin{bmatrix} a\,b \\ 1 \end{bmatrix} \otimes x'\,\underset{w}{\overset{y'}{\boxslash}}\,z$$

(applying an associativity transformation of lower level to the second component)

$$\cong \bigoplus_{u\neq 1} \begin{bmatrix} a\,b \\ u \end{bmatrix} \otimes \begin{bmatrix} x'uy'\ z \\ w \end{bmatrix} \quad + \quad \begin{bmatrix} a\,b \\ 1 \end{bmatrix} \otimes x'\,\underset{w}{\overset{y'}{\boxbslash}}\,z$$

$$= \sum_{\substack{u\neq 1 \\ v\neq 1}} \delta_{x'uy''vz',w} \begin{bmatrix} a\,b \\ u \end{bmatrix} \otimes \begin{bmatrix} b'\ c \\ v \end{bmatrix} \quad + \quad \sum_{u\neq 1} \begin{bmatrix} a\,b \\ u \end{bmatrix} \otimes \begin{bmatrix} b'\ c \\ 1 \end{bmatrix} \otimes \begin{bmatrix} x'uy''\ z' \\ w \end{bmatrix}$$

$$+ \sum_{v\neq 1} \begin{bmatrix} a\,b \\ 1 \end{bmatrix} \otimes \begin{bmatrix} b'\ c \\ v \end{bmatrix} \otimes \begin{bmatrix} x'\ y''vz' \\ w \end{bmatrix} \quad + \quad \begin{bmatrix} a\,b \\ 1 \end{bmatrix} \otimes \begin{bmatrix} b'\ c \\ 1 \end{bmatrix} \otimes x'\,\underset{w}{\overset{y''}{\boxbslash}}\,z'$$

(using the coherence for one-vertex interaction)

$$\cong \sum_{\substack{u\neq 1\\ v\neq 1}} \delta_{x'uy''vz',w} \begin{bmatrix} a\,b \\ u \end{bmatrix} \otimes \begin{bmatrix} b'\ c \\ v \end{bmatrix} \quad + \quad \sum_{u\neq 1} \begin{bmatrix} a\,b \\ u \end{bmatrix} \otimes \begin{bmatrix} b'\ c \\ 1 \end{bmatrix} \otimes \mathrm{Hom}(x'uy''z', w)$$

$$+ \sum_{v\neq 1} \begin{bmatrix} a\,b \\ 1 \end{bmatrix} \otimes \begin{bmatrix} b'\ c \\ v \end{bmatrix} \otimes \mathrm{Hom}(x'y''vz', w) \quad + \quad \begin{bmatrix} a\,b \\ 1 \end{bmatrix} \otimes \begin{bmatrix} b'\ c \\ 1 \end{bmatrix} \otimes$$ 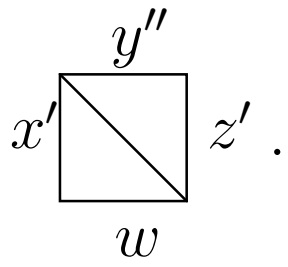

The last expression is simply denoted by

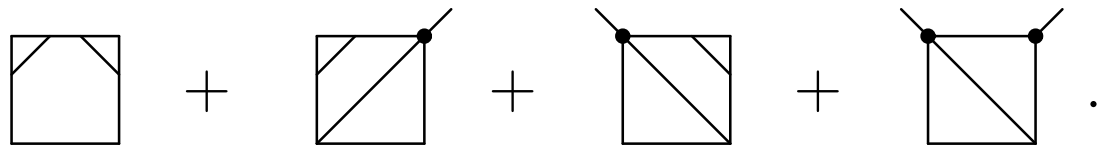

Similarly we have

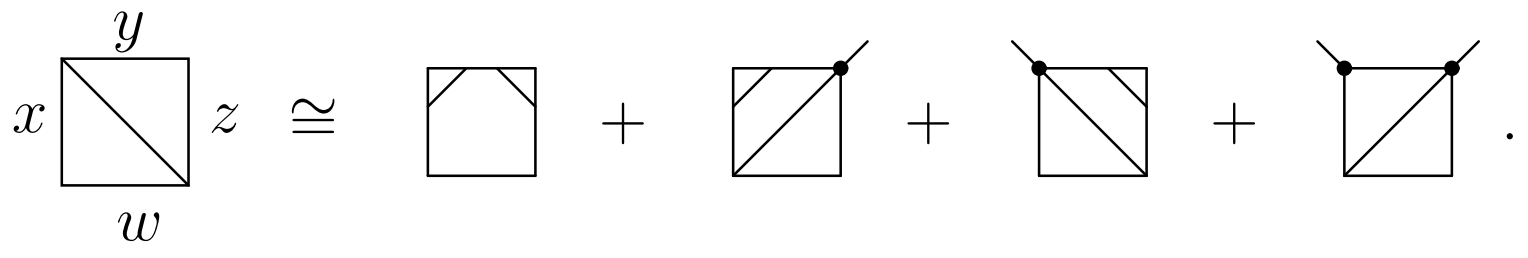

Now the isomorphism 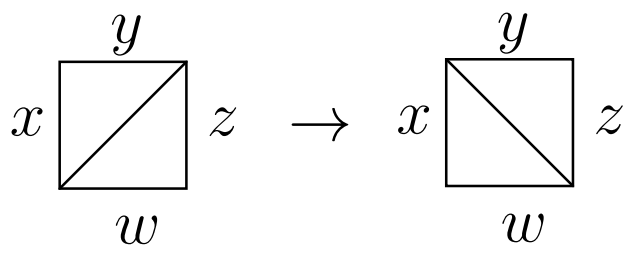 is defined by applying the associativity transformation (of lower level) to 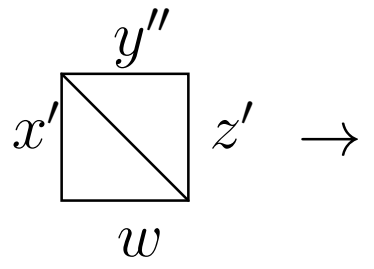

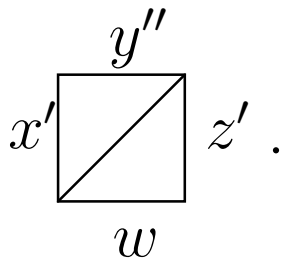

## 2. Pentagonal Relations

We shall here prove the pentagonal relations for the associativity transformations introduced in the previous section.

Consider a pentagon labeled by elements in $S * T$. If there appears the unit element in the labeling of edges other than the bottom edge,

then the coherence is satisfied as can be seen below:

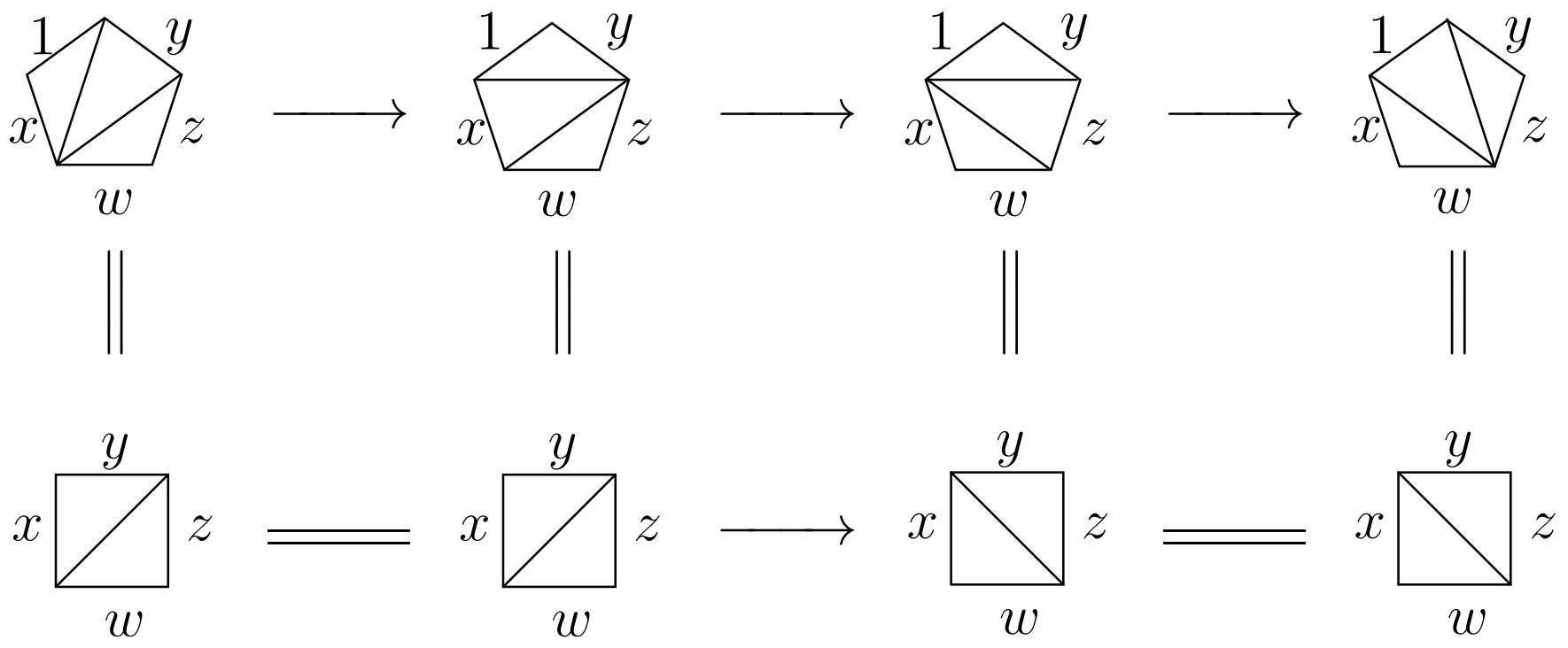

and

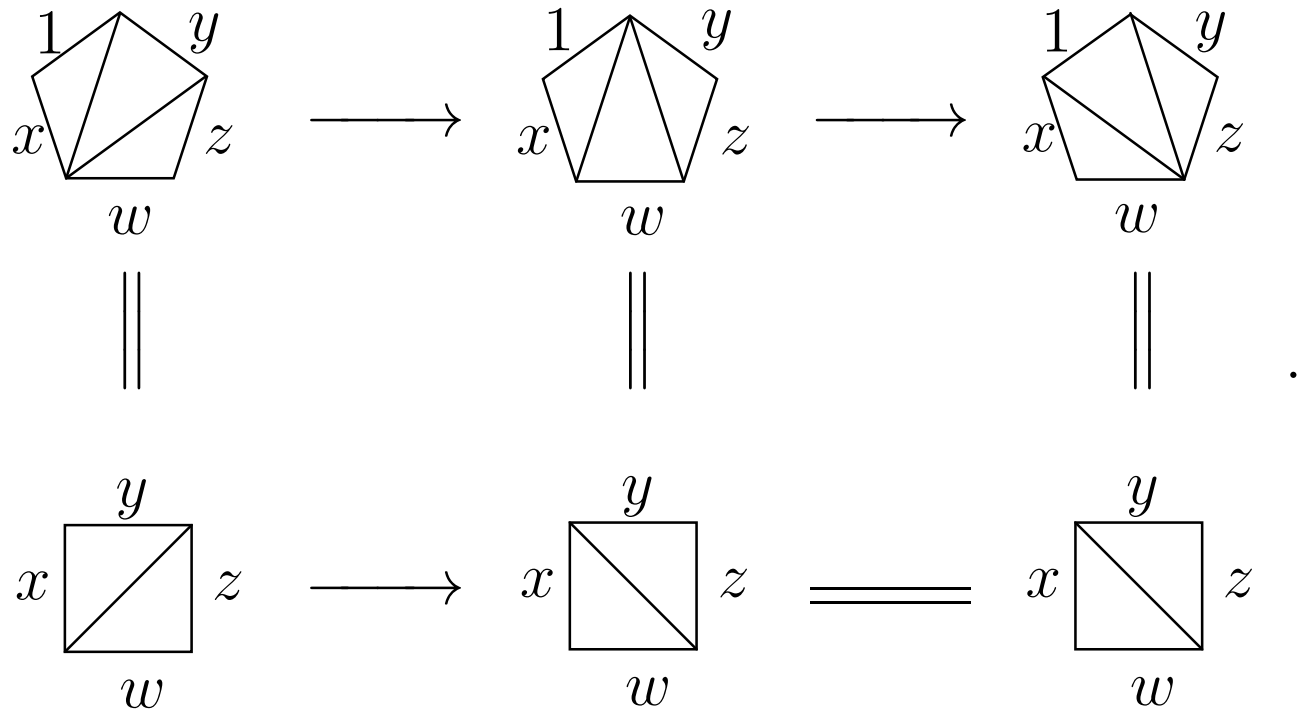

Thus we may assume that the unit 1 does not appear in the labeling. In that case, we add corners at intermediate points of each label and enlarge the pentagon into a polygon (with the number of edges except for the bottom given by the total length of labeling) which has the labeling by elements in $S^\circ \sqcup T^\circ$ with at most three vertices of interaction.

We now prove the coherence for labeled polygons of this type by an induction on the length of labeling (the label of the bottom edge being not counted for the length).

If all the labels belong to the same group ($S^\circ$ or $T^\circ$), the problem is reduced to the coherence in $\mathcal{C}$ or $\mathcal{D}$ respectively. Therefore we need to consider the case that one of vertices is inactive. Now choose one of such vertices (denoted by $*$) and use it as the reference point in the proof of coherence in Appendix A. We are then faced with the coherence for subpentagons containing $*$ as a vertex.

The labeling of such subpentagons (which really contribute to the polygonal vector space) should contain at most two interactive vertices unless it has a smaller length. In fact, if there exists a labeling of the same length, the subpentagon must contain the bottom edge of the polygon because contractions do not increase the length of labels (cf. Lemma 1.1). If it further contains three interactive points (denoted $\bullet$ in Fig. 3), then one of them should be in the middle of the arc cut by an edge of the subpentagon (the vertex $*$ being inactive), whence non-trivial contractions take place for labels of smaller length at that edge, contradicting with the length assumption.

As we can apply the induction hypothesis for labelings of smaller length, there remains to consider the labelings of subpentagons which have at most two interactive points located at vertives other than $*$ (the left and right shoulder vertices of pentagon).

If it has just one interactive point, we can apply the one-vertex coherence Lemma 1.3 and hence we are lead to the case of two interactive points located at left and right shoulder vertices of pentagon.

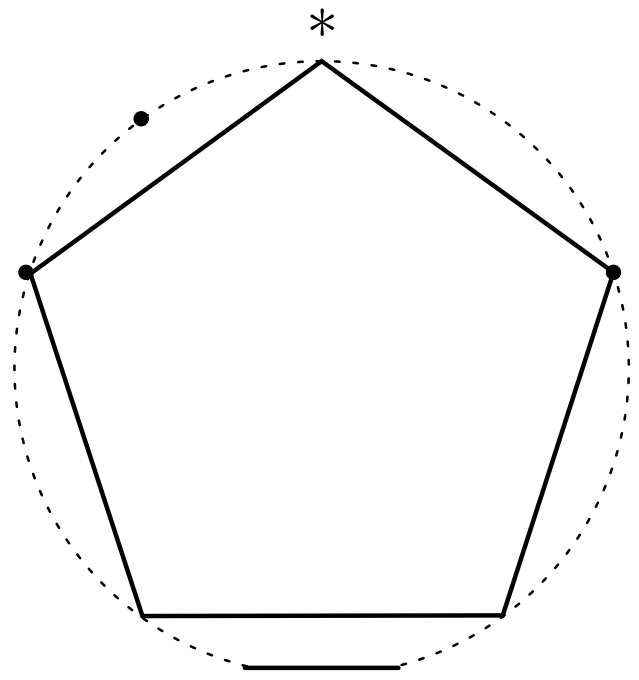

Figure 3.

Consider the labeling decsribed by 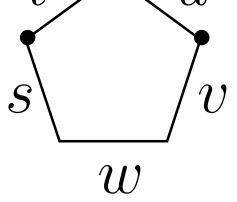

with $s|t$, $t||u$ and $u|v$ ($\bullet$ indicates interactive points). Write $s = s'a$, $t = bt'$, $u = u'x$ and

$v = yv'$ with $|a| = |b| = |x| = |y| = 1$. Then we have

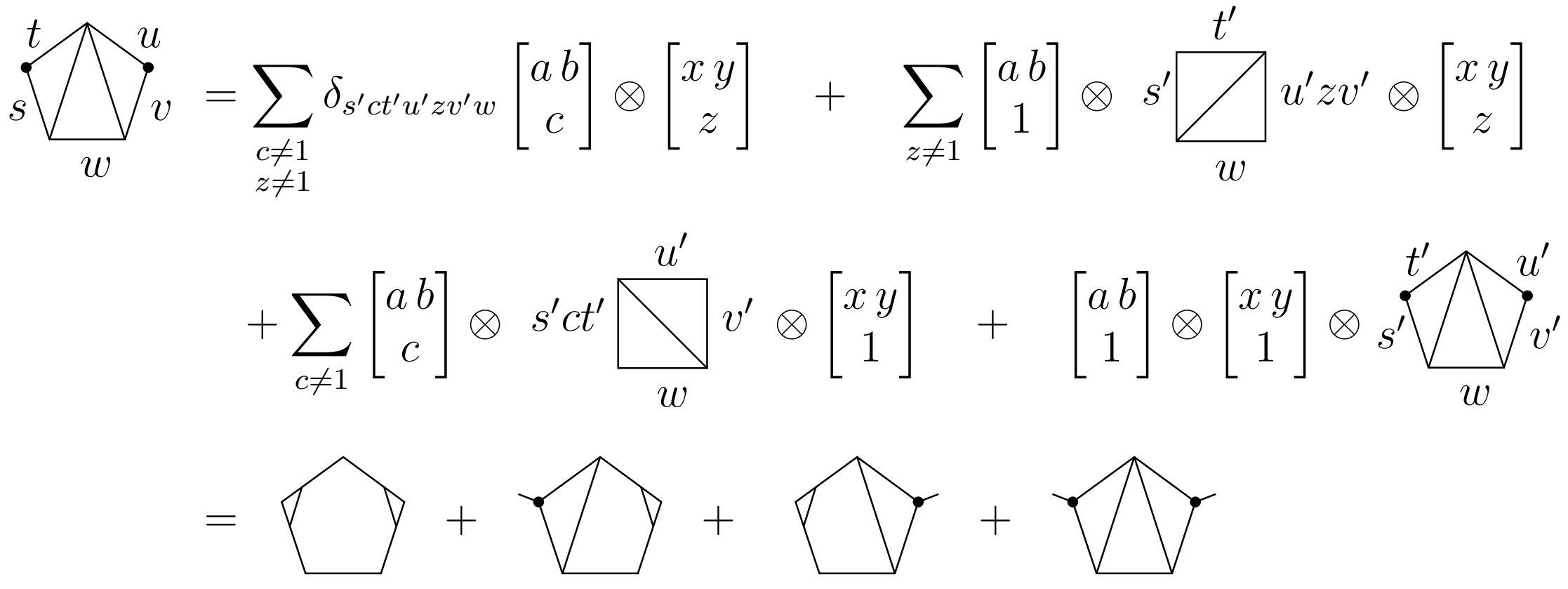

Next we define the isomorphism

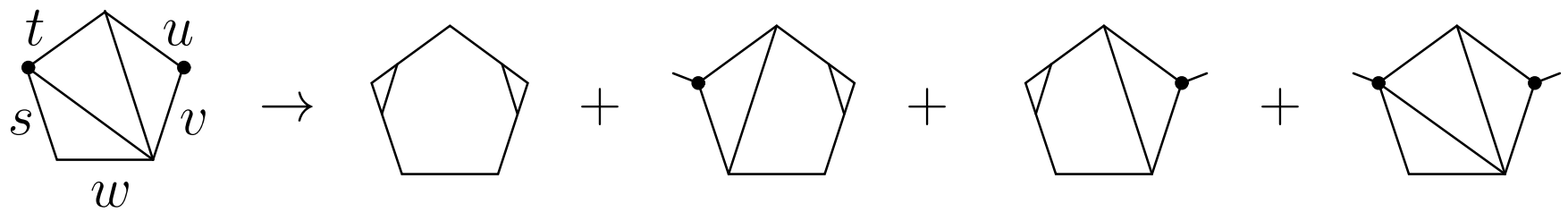

by the commutativity of the diagram

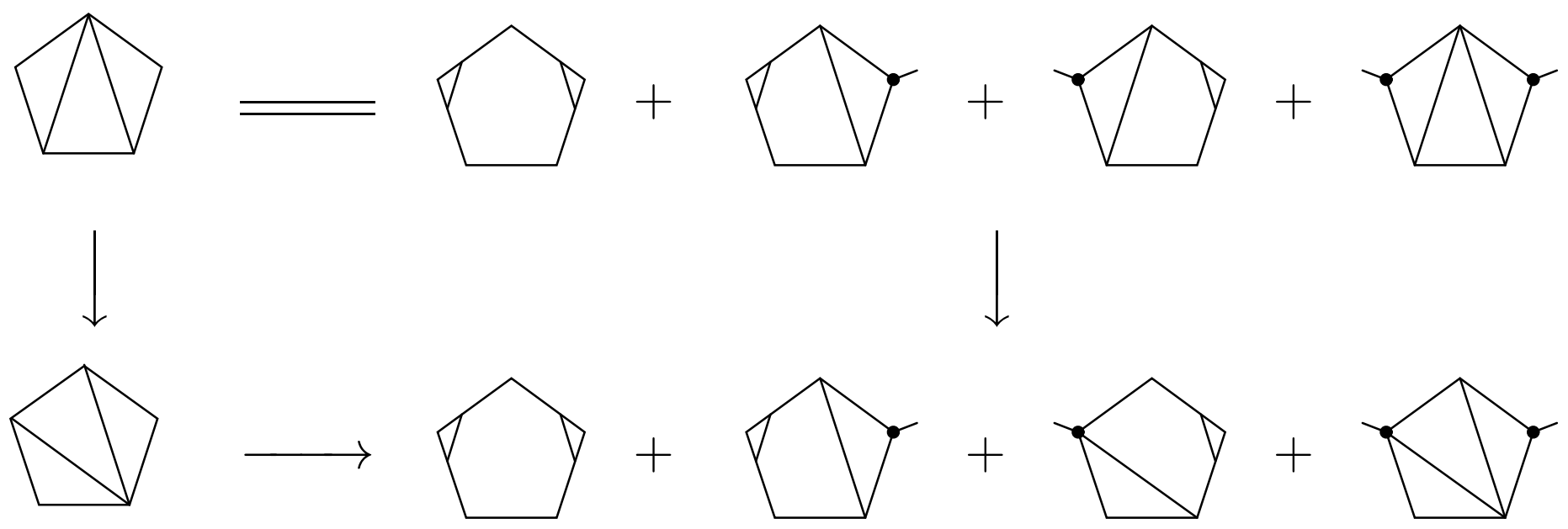

and then introduce the isomorphism

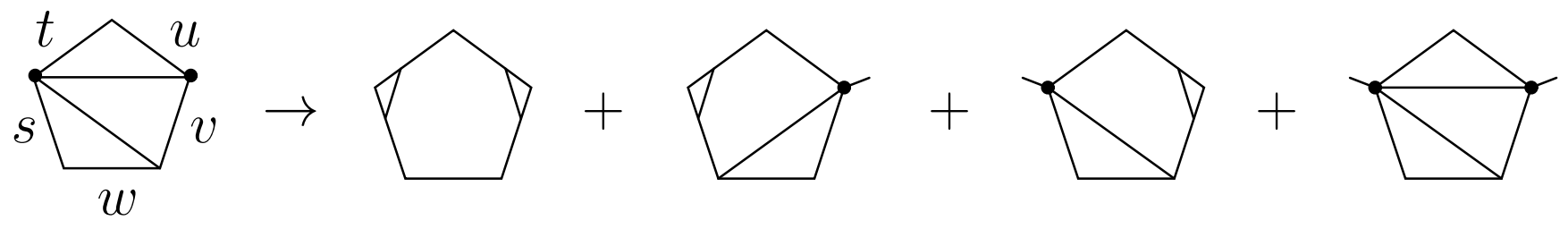

by

(1)

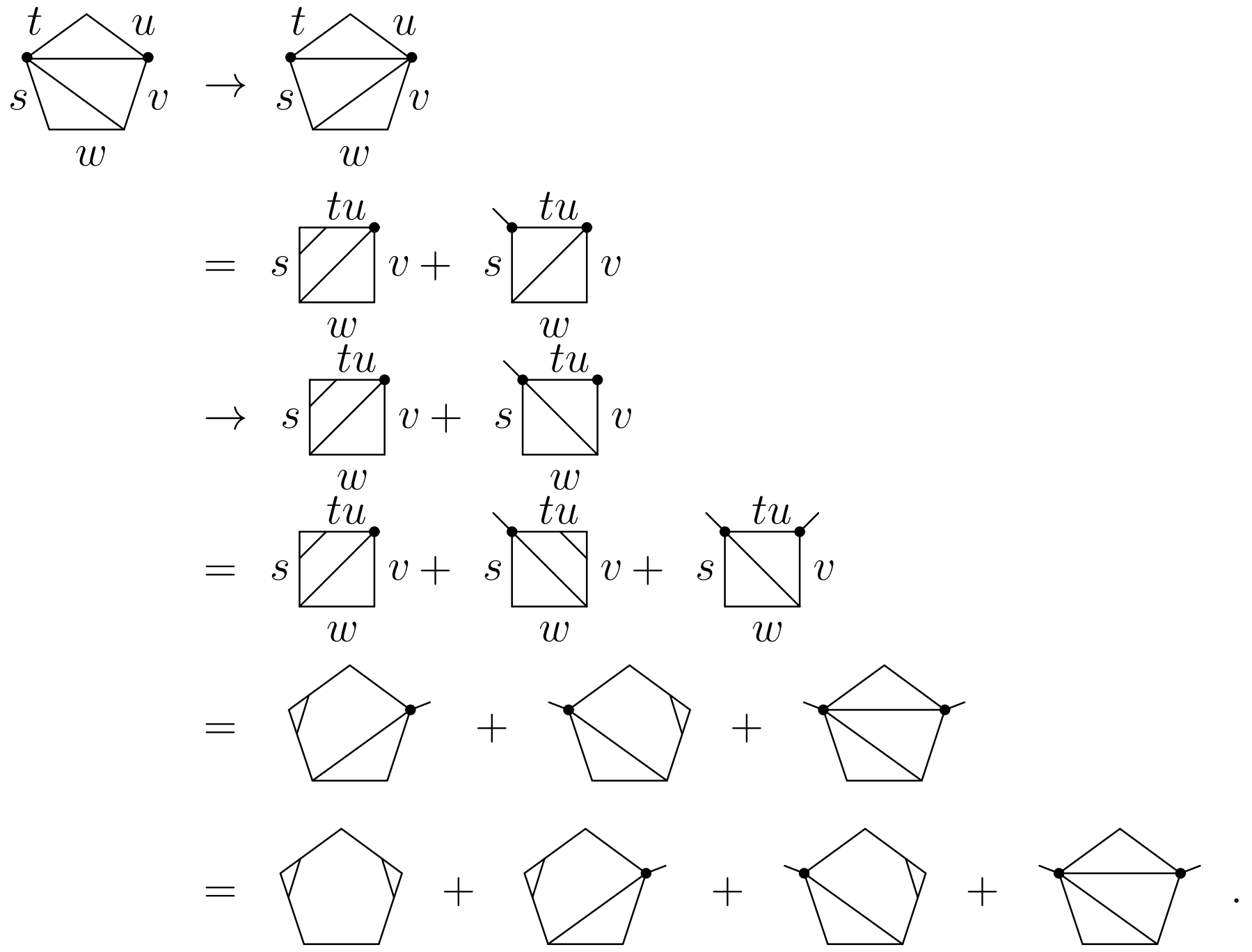

Note here that the associativity transformation 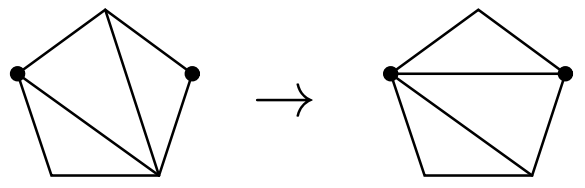
is the identity when they are expressed by triangulated vector spaces (Lemma 1.2).

Now we claim the commutativity of the diagram

(2)

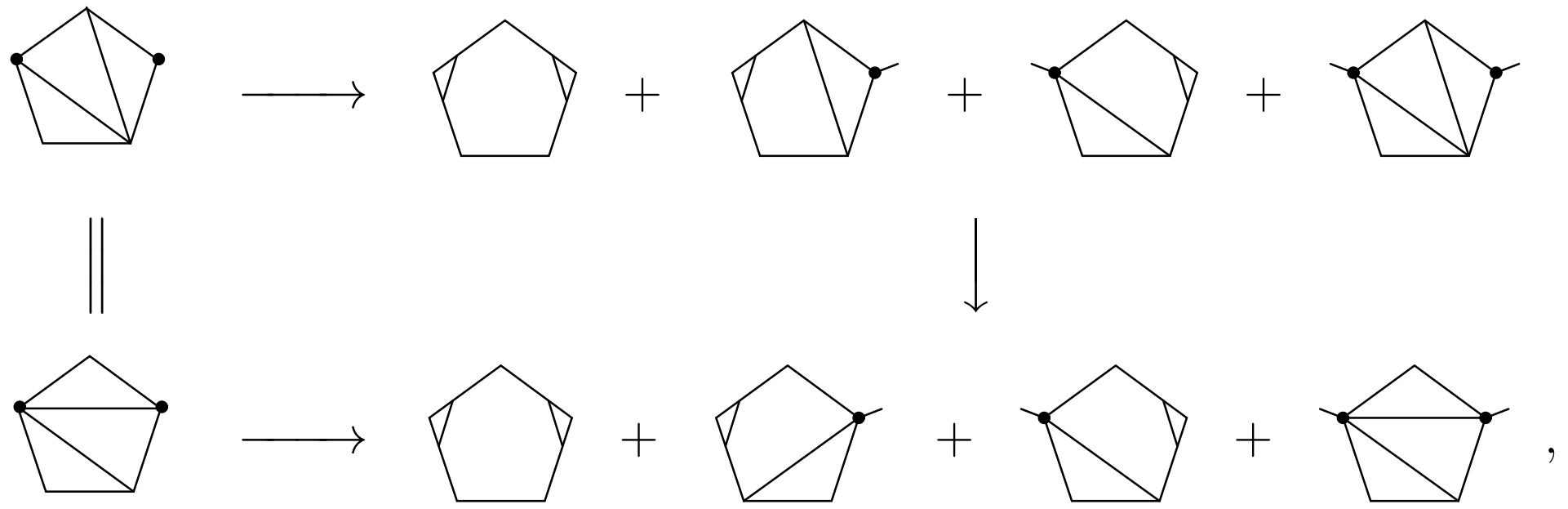

where the right vertical arrow represents associativity transformations of lower level.

We shall prove this by an induction on the length $|u|$. From the inductive formula for 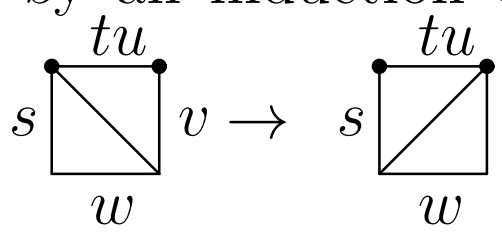

, we have

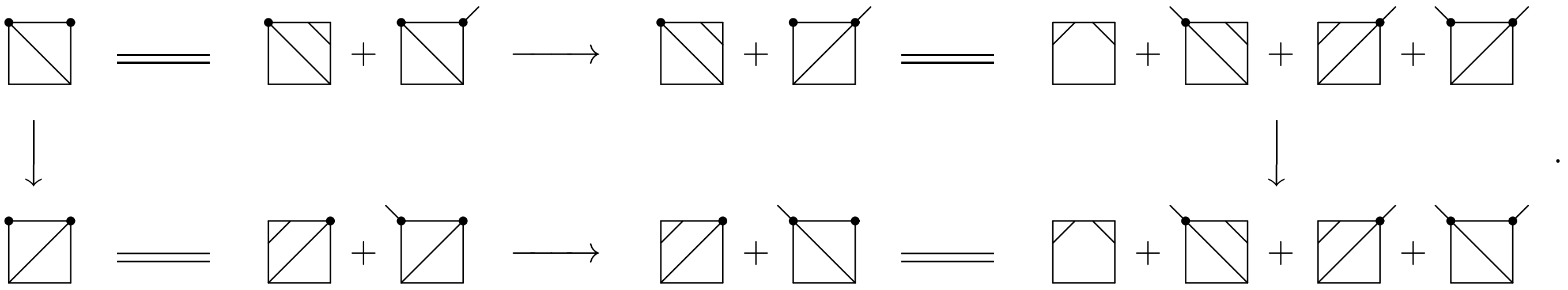

If we combine this with the definition of the 4-component decomposition of 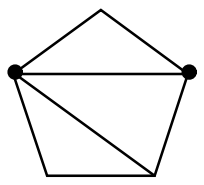 by using the identities

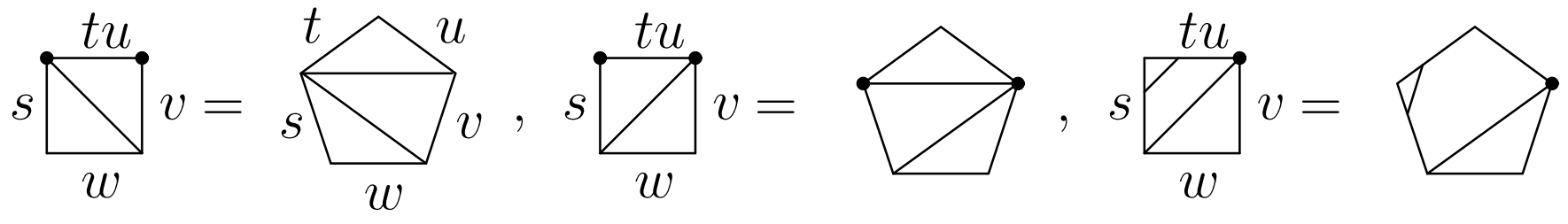

and so on, the isomorphism

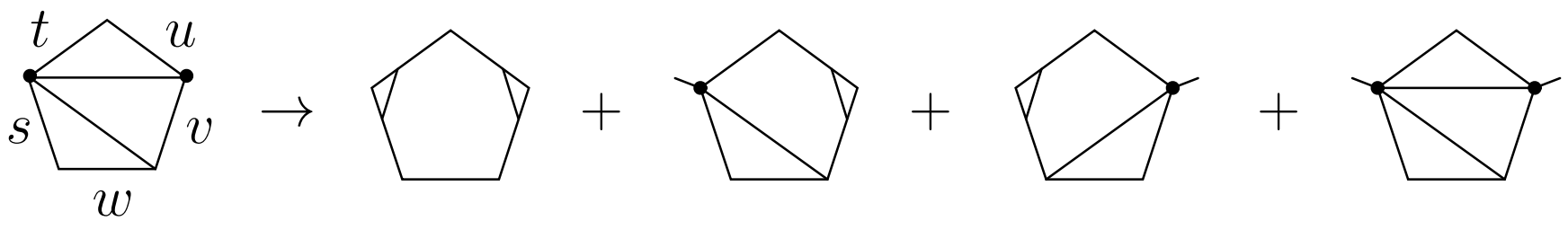

turns out to be given by

(3)

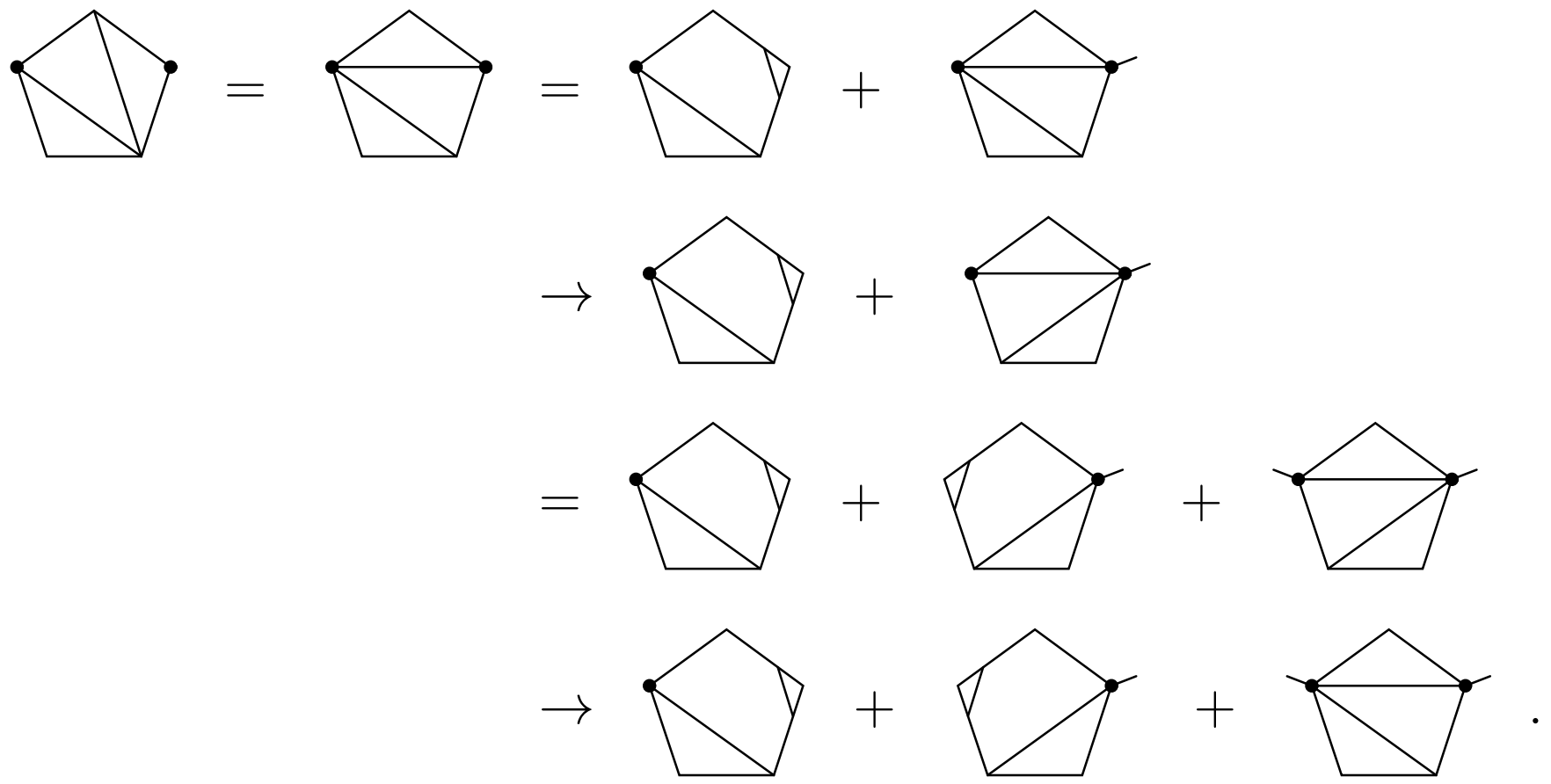

Viewing the definition of the 4-component decomposition of 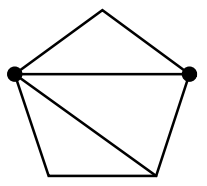

(see (1)), the changes on 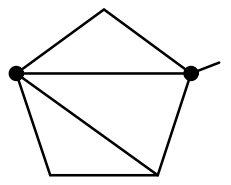 in the above transformation are exactly the ones given by (1).

If $|u| = 1$, the transformation (3) is reduced to

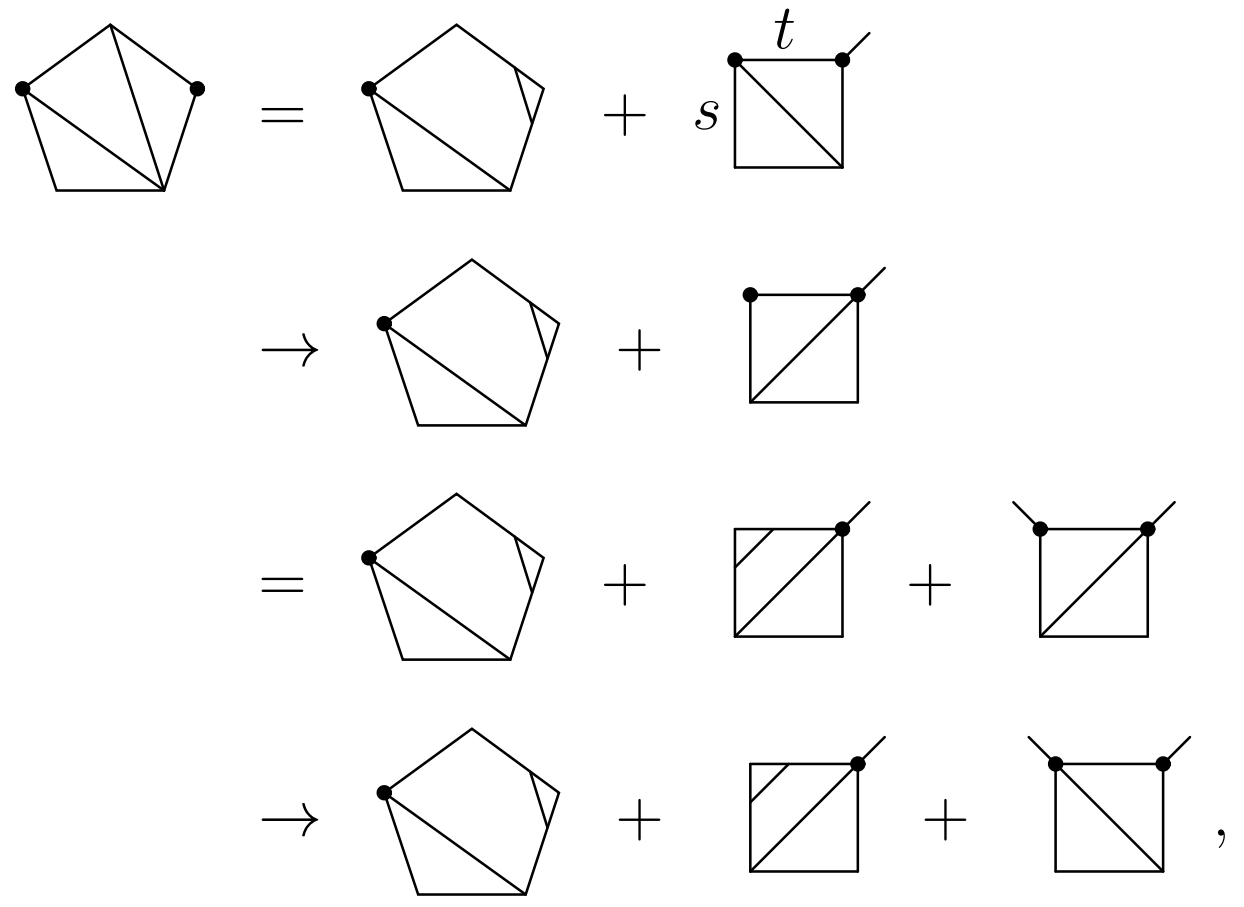

which is nothing but the composite isomorphism

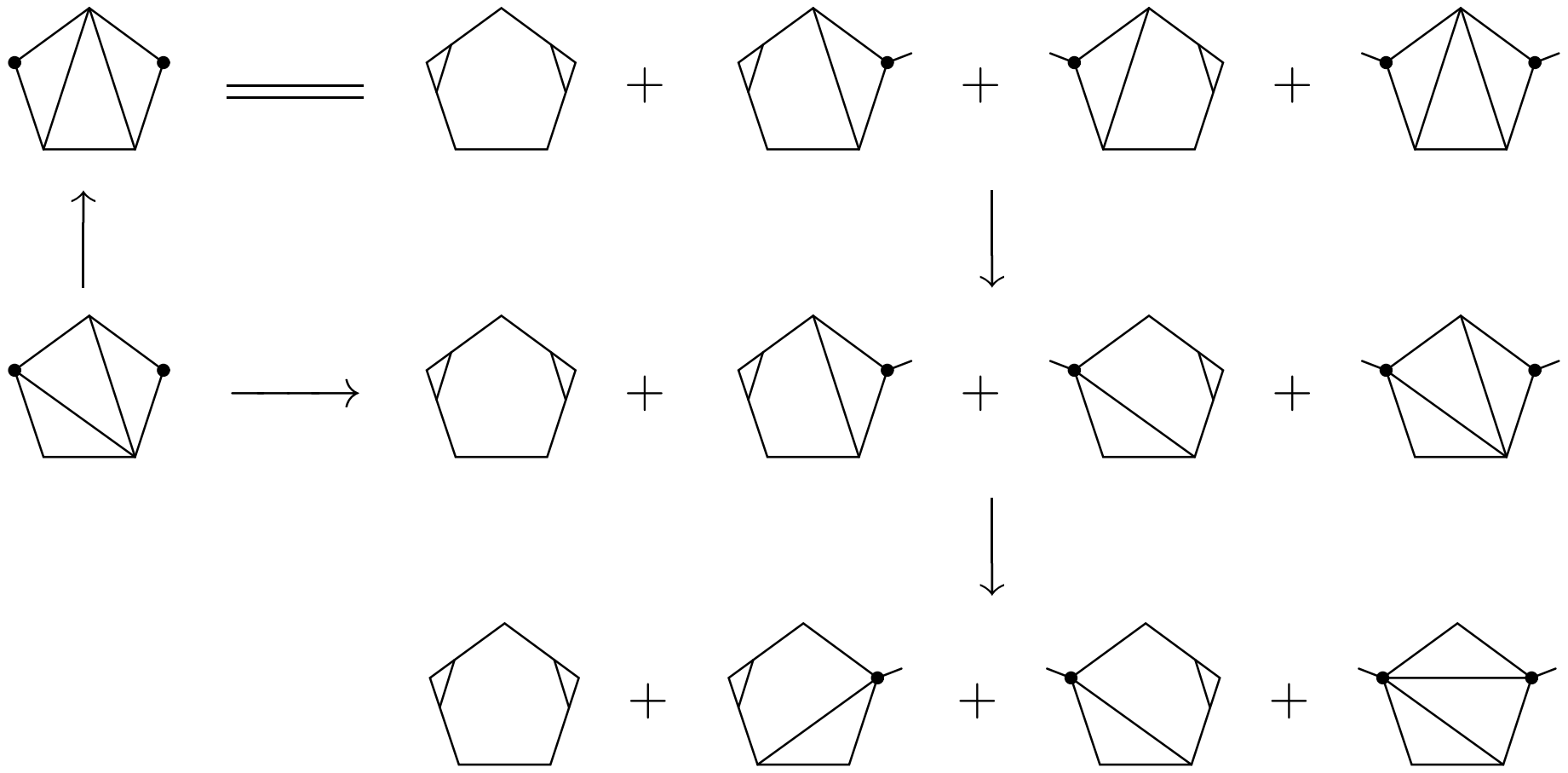

for $|u| = 1$, showing the commutativity of the diagram in question.

For $|u| \geq 2$, assume the induction hypothesis that the diagram commutes for lower levels to replace the operation of 4-component decomposition on 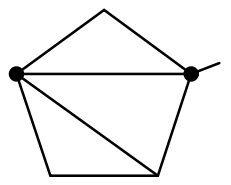 inside (3) by the above composite isomorphism. As a consequence, we obtain the following expression for the 4-component decomposition of 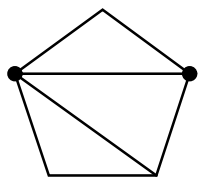 .

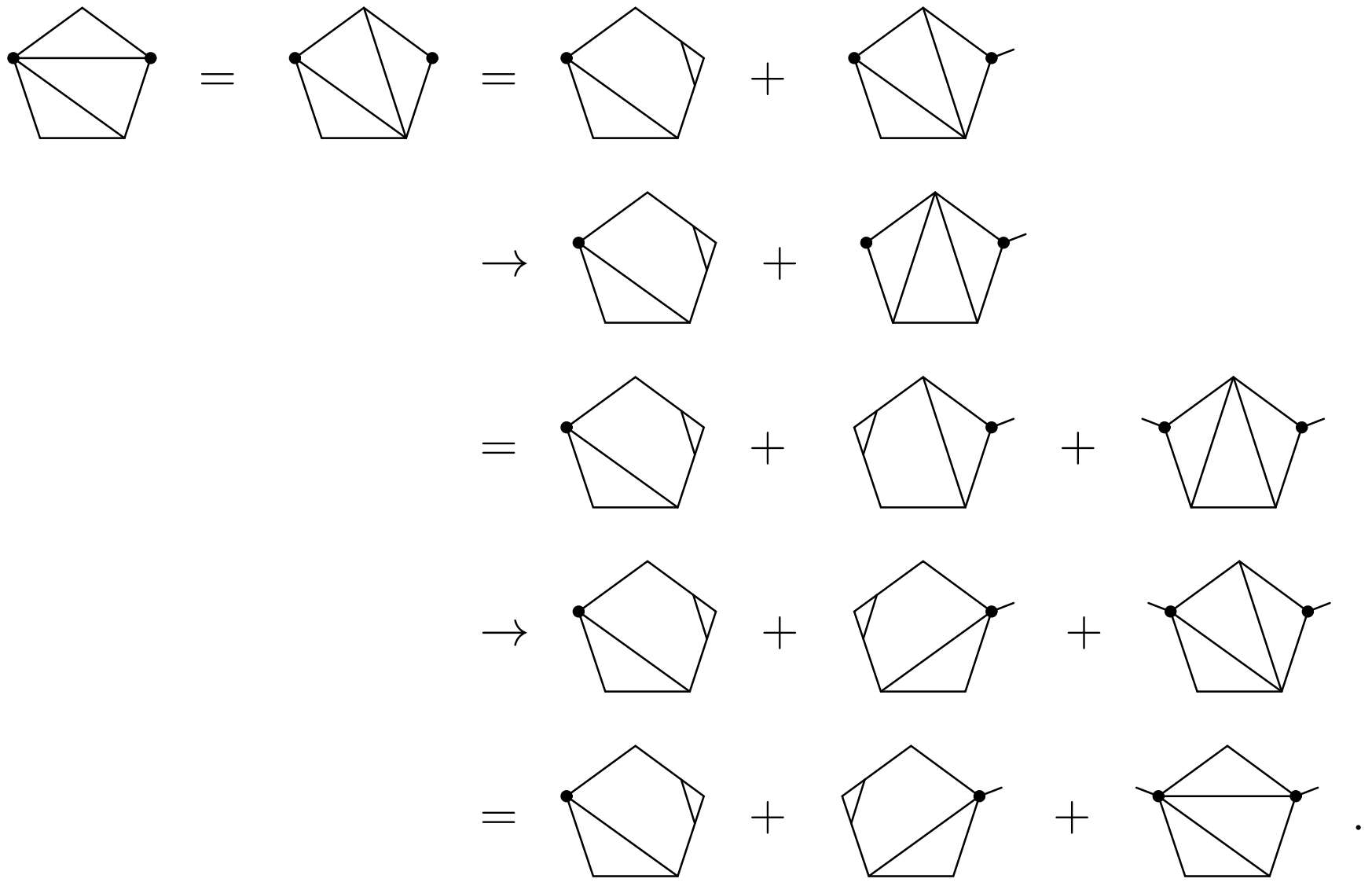

From the inductive relation of associativity transformations of one-point interaction, we have the commutativity and so on, whence the above isomorphism is the one given by the path directed to right and then down in (2), proving the commutativity of the diagram (2).

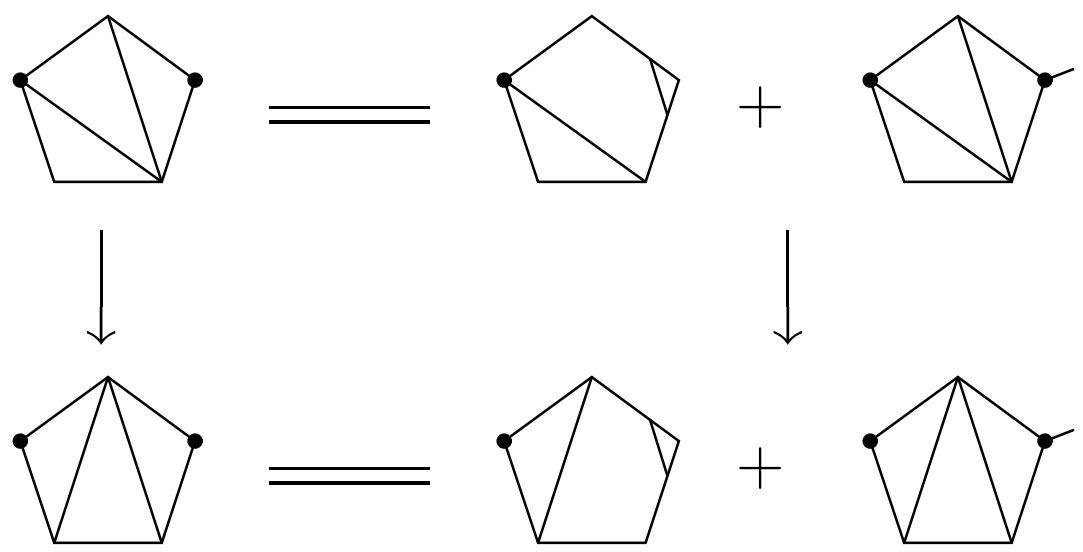

Finally we introduce the 4-component decomposition for 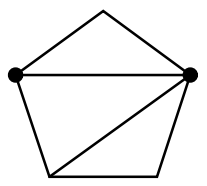

by

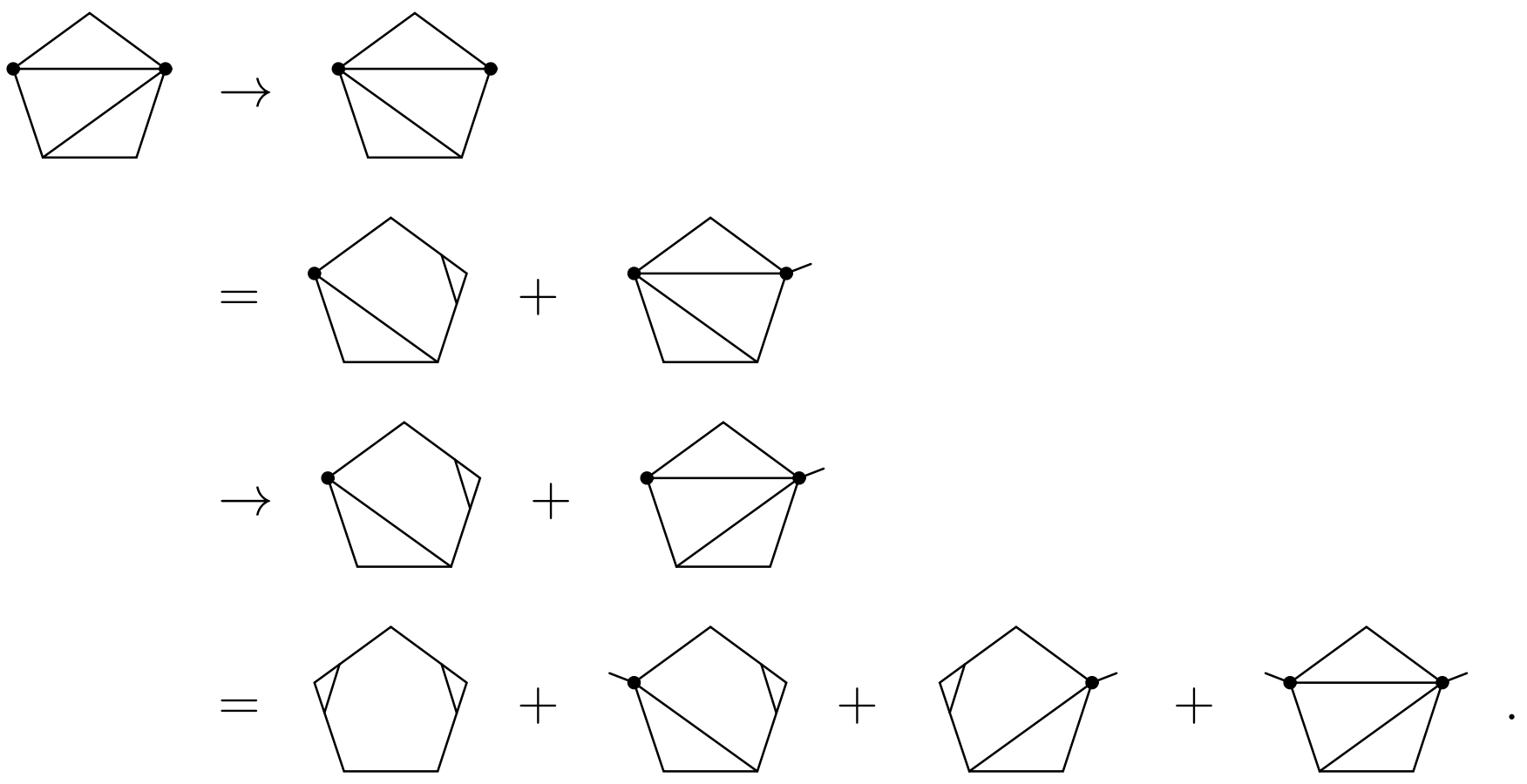

The commutativity of the diagram

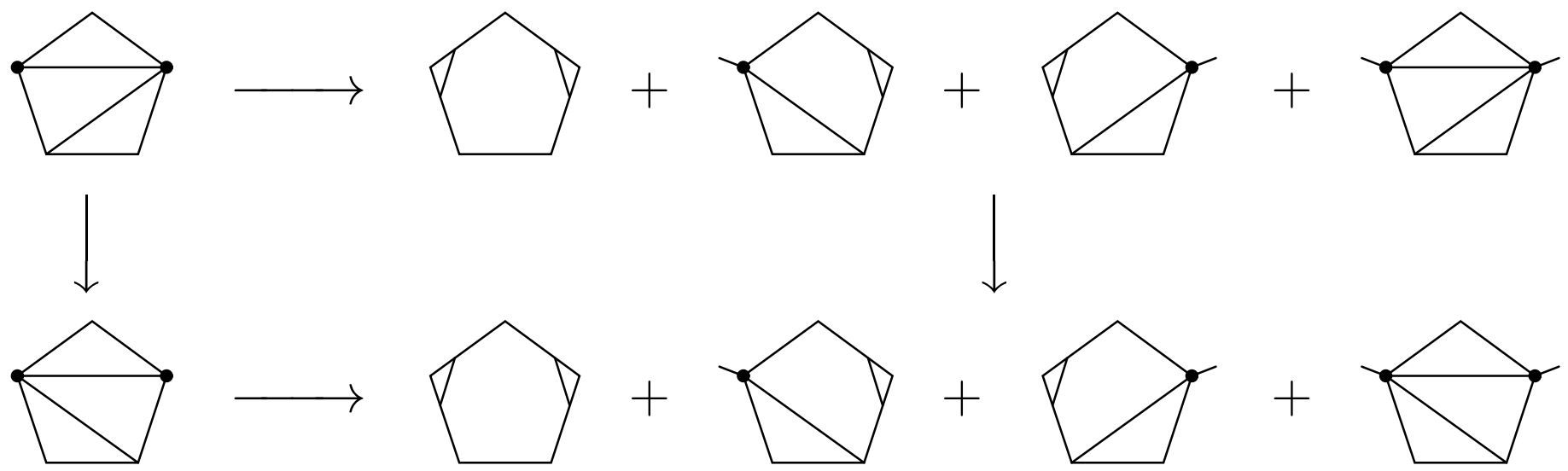

is then reduced to that of

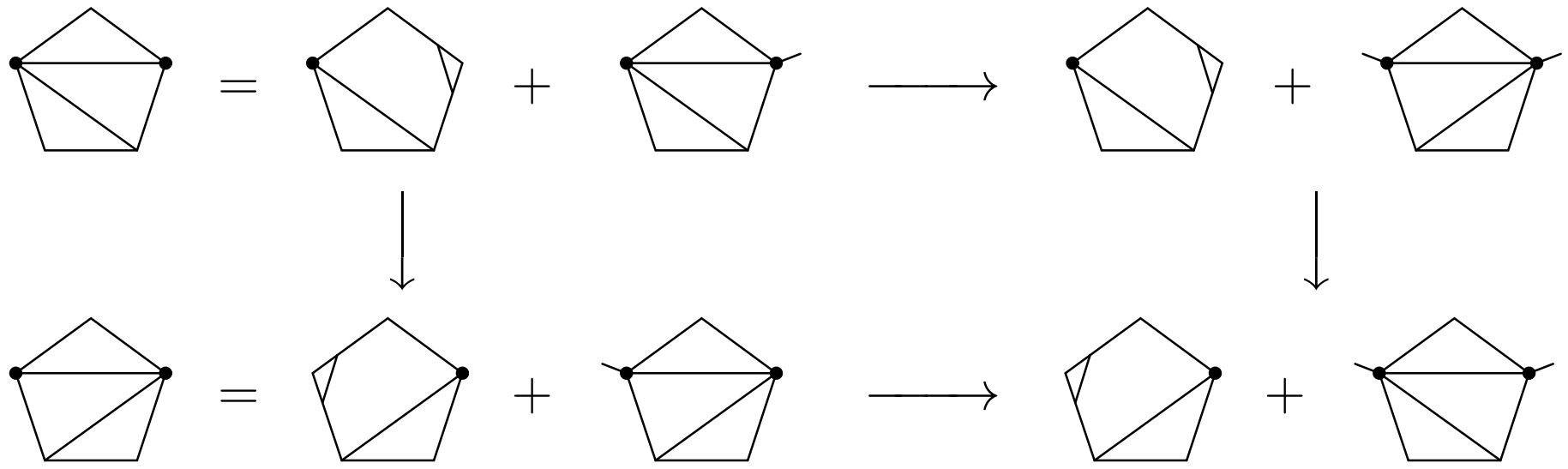

which is nothing but the induction formula of associativity transformation for $s$ 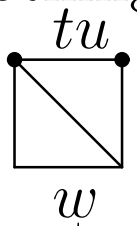

$v \to s$ 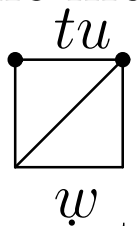

$v$.

By symmetry, we can introduce isomorphisms

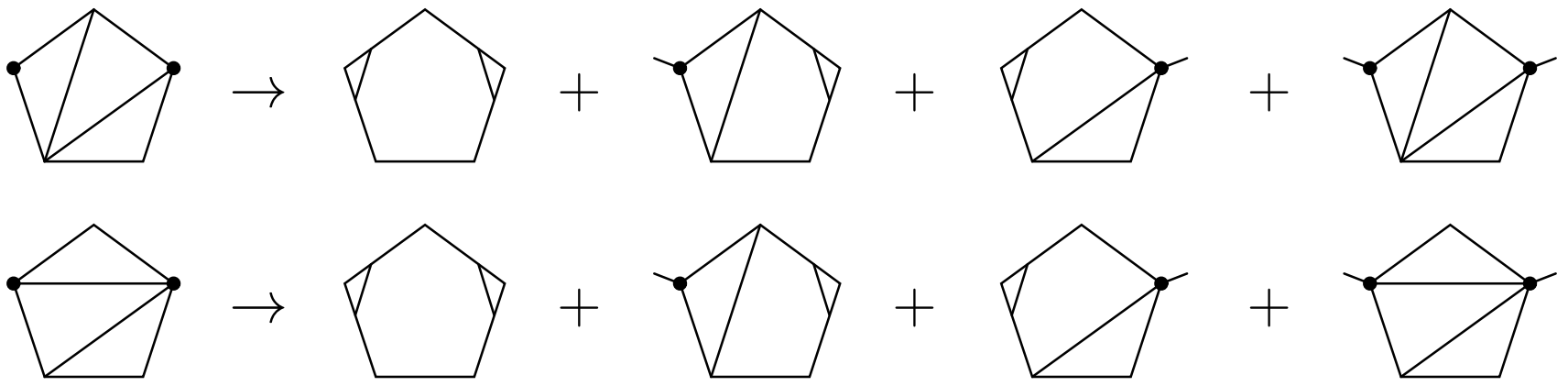

so that they form the corresponding commutative diagrams.

With all these isomorphisms in hand, we can now prove the pentagonal relation for the associativity transformations: Putting these commutative diagrams of 4-component decompositions together, we obtain the commutative diagram Fig. 4, where the outcircuit is the pentagonal relation for 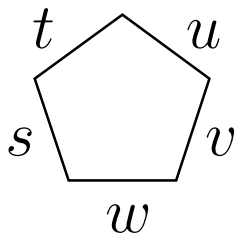

, spokes represent 4-component decompositions and the inner pentagon is given by associativity transformations on

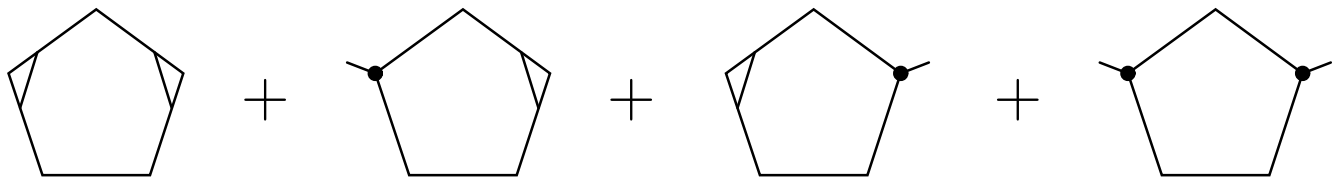
.

The operations are the identities on 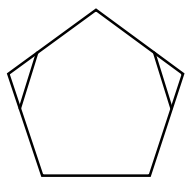 , transformations of one-point interaction on 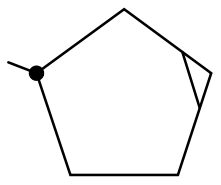 , 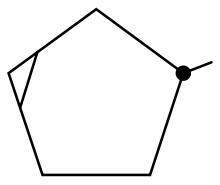 and associativity transformations in the pentagonal relation for 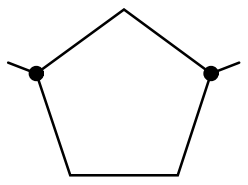 . By Lemma 1.3, the

circuits for 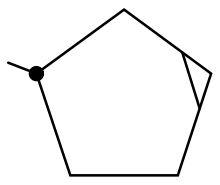 and 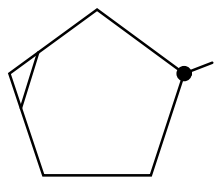 commute, whereas the pentagonal relation for 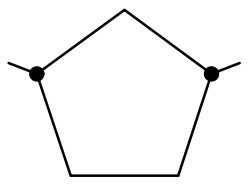 is satisfied by the induction hypothesis on the length $|s| + |t| + |u| + |v|$.

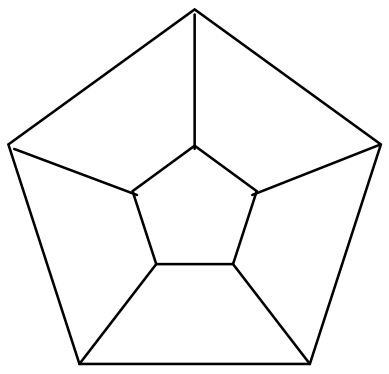

FIGURE 4.

## 3. FREE PRODUCTS

Being preprated in the previous section, we here define the free products of semisimple tensor categories, say $\mathcal{C}$ and $\mathcal{D}$ with the spectrum sets denoted by $S$ and $T$ respectively.

We then have the fusion algebra $\mathbb{Z}[S * T]$, triangular vector spaces $\left\{ \begin{bmatrix} x\,y \\ z \end{bmatrix} \right\}_{x,y,z \in S*T}$ and a family $\{a_{x,y,z,w}\}_{x,y,z,w \in S*T}$ of isomorphisms satisfying the pentagonal relation with

$$a_{x,y,z,w} : \bigoplus_{u \in S*T} \begin{bmatrix} x\,y \\ u \end{bmatrix} \otimes \begin{bmatrix} u\,z \\ w \end{bmatrix} \quad \to \quad \bigoplus_{v \in S*T} \begin{bmatrix} y\,z \\ v \end{bmatrix} \otimes \begin{bmatrix} x\,v \\ w \end{bmatrix}.$$

Recall that the triangular vector spaces are such that

$$\begin{bmatrix} 1\,x \\ y \end{bmatrix} = \mathbb{C}\delta_{x,y} = \begin{bmatrix} x\,1 \\ y \end{bmatrix}$$

with the isomorphism $a_{x,1,y,z}$ satisfying

$$a_{x,1,y,z}(1_x \otimes \sigma) = 1_y \otimes \sigma \quad \text{for } \sigma \in \begin{bmatrix} x\,y \\ z \end{bmatrix},$$

where $1_x$ and $1_y$ denote the unit element of $\mathbb{C}$ in triangular vector spaces

$$\begin{bmatrix} x\,1 \\ x \end{bmatrix}, \qquad \begin{bmatrix} 1\,y \\ y \end{bmatrix}$$

respectively.

As described in [16], these in fact (re)construct the tensor category $\mathcal{C} * \mathcal{D}$, the free product of $\mathcal{C}$ and $\mathcal{D}$: Objects of $\mathcal{C} * \mathcal{D}$ are families of finite-dimensional complex vector spaces indexed by the set $S * T$, say $X = \{X(u)\}_{u \in S*T}$, which are essentially finite, i.e., $X(u) = \{0\}$ for $u \notin F$ with $F$ a finite subset of $S * T$. The hom-sets are defined by

$$\mathrm{Hom}(X, Y) = \bigoplus_{u \in S*T} \mathrm{Hom}(X(u), Y(u)).$$

The tensor product is given by

$$(X \otimes Y)(u) = \bigoplus_{x,y \in S*T} \begin{bmatrix} x\, y \\ u \end{bmatrix}^* \otimes X(x) \otimes Y(y)$$

with the unit object $I$ by

$$I(u) = \mathbb{C}\delta_{u,1}.$$

For morphisms $f : X \to X'$ and $g : Y \to Y'$, we set

$$(f \otimes g)(u) = \bigoplus_{x,y \in S*T} \left( \begin{bmatrix} x\, y \\ u \end{bmatrix}^* \otimes X(x) \otimes Y(y) \xrightarrow{1 \otimes f(x) \otimes g(y)} \begin{bmatrix} x\, y \\ u \end{bmatrix}^* \otimes X'(x) \otimes Y'(y) \right).$$

The left and right unit constaints are then defined by the obvious identities

$$(l_X)(u) : (I \otimes X)(u) = \mathbb{C} \otimes X(u) \to X(u),$$
$$(r_X)(u) : (X \otimes I)(u) = X(u) \otimes \mathbb{C} \to X(u)$$

and the associativity transformations by

$$\bigoplus_{x,y,z \in S*T} \left( \begin{pmatrix} & y & \\ x & \diagup & z \\ & w & \end{pmatrix}^* \otimes X(x) \otimes Y(y) \otimes Z(z) \xrightarrow{a_{x,y,z,w} \otimes 1} \begin{pmatrix} & y & \\ x & \diagdown & z \\ & w & \end{pmatrix}^* \otimes X(x) \otimes Y(y) \otimes Z(z) \right.$$

with

$$((X \otimes Y) \otimes Z)(w) = \bigoplus_{x,y,z \in S*T} \begin{pmatrix} & y & \\ x & \diagup & z \\ & w & \end{pmatrix}^* \otimes X(x) \otimes Y(y) \otimes Z(z),$$

$$(X \otimes (Y \otimes Z))(w) = \bigoplus_{x,y,z \in S*T} \begin{pmatrix} & y & \\ x & \diagdown & z \\ & w & \end{pmatrix}^* \otimes X(x) \otimes Y(y) \otimes Z(z).$$

The tensor categories $\mathcal{C}$ and $\mathcal{D}$ are naturally identified with tensor subcategories of $\mathcal{C} * \mathcal{D}$: For an object $X$ in $\mathcal{C}$, the associated object $\{X(u)\}_{u \in S*T}$ in $\mathcal{C} * \mathcal{D}$ is defined by

$$X(u) = \begin{cases} \mathrm{Hom}(u, X) & \text{if } u \in S, \\ \{0\} & \text{otherwise.} \end{cases}$$

For a morphism $f : X \to X'$ in $\mathcal{C}$, we set

$$f(u) : X(u) \ni \xi \mapsto f \circ \xi \in X'(u).$$

The tensor category $\mathcal{C} * \mathcal{D}$ is then monoidally equivalent to the one generated by its subcategories $\mathcal{C}$ and $\mathcal{D}$.

When $\mathcal{C}$ and $\mathcal{D}$ are rigid, their free product is rigid because tensor products of rigid objects are rigid as well.

If $\mathcal{C}$ and $\mathcal{D}$ are C*-tensor categories, the free product $\mathcal{C}*\mathcal{D}$ is naturally a C*-tensor category: the triangular vector spaces are Hilbert spaces as direct sums of finitely many tensor products of Hilbert spaces. The associativity transformations are then unitary because the recursive definition gives the same type of structure: direct sums of tensor products of associativity transformations in $\mathcal{C}$ or $\mathcal{D}$. The left and right unit constraint vectors are clearly normalized.

When $\mathcal{C}$ and $\mathcal{D}$ are rigid, the free product C*-tensor category $\mathcal{C} * \mathcal{D}$ is rigid and hence it admits a positive Frobenius duality (see [23] for an explicit description in terms of polygonal presentation).

As a final remark, we record here that our construction is obviously extented to the free product of an arbitrary family of semisimple tensor categories: given a family $\{\mathcal{C}_\alpha\}_{\alpha \in A}$ of semisimple tensor categories, its free product $*_{\alpha \in A} \mathcal{C}_\alpha$ is defined in such a way that

$$*_{\alpha \in A} \mathcal{C}_\alpha = *_{j \in J} \left( *_{\alpha_j \in A_j} \mathcal{C}_{\alpha_j} \right),$$

where $\{A_j\}_{j \in J}$ denotes a division of the index set $A$.

*Remark* . If tensor categories $\mathcal{C}$ and $\mathcal{D}$ are realized in the tesnor category $\mathcal{V}$ of finite-dimensional vector spaces (i.e., they are the Tannaka duals of Hopf algebras), then we have an obvious realization of the free product $\mathcal{C}*\mathcal{D}$ in $\mathcal{V}$ and hence $\mathcal{C}*\mathcal{D}$ itself is the Tannaka dual of a Hopf algebra. This fact, however, can be deduced more easily by dealing with free products of Hopf algebras and their corepresentations as worked out by S. Wang in [17].

## 4. Bisch-Jones' Planar Algebras

In this section, we shall apply our results to identify the planar algebras of Bisch and Jones.

Let $S = \{s_n; n \geq 0\}$ be the fusion rule set of (the Tannaka dual of) $SU(2)$:

$$s_m s_n = \sum_{k=0}^{\min\{m,n\}} s_{|m-n|+2k}, \qquad s_n^* = s_n.$$

According to Kazhdan and Wenzl ([12]), rigid semisimple tensor categories (over the field $\mathbb{C}$) having the fusion rule set $S$ is classified in the following way: Let $X$ be an object representing $s_1 \in S$. Then, by the fusion rule, we have

$$X \otimes X \cong I \oplus X'$$

with $X'$ representing $s_2$. Let $e \in \mathrm{End}(X \otimes X)$ be the idempotent to the $I$-component and set $e_1 = e \otimes 1_X$, $e_2 = 1_X \otimes e$. Then we have

$$e_1 e_2 e_1 = (q + q^{-1} + 2)^{-1} e_1, \quad e_2 e_1 e_2 = (q + q^{-1} + 2)^{-1} e_2,$$

with $q \in \mathbb{C}$ satisfying $1 + q + q^2 + \cdots + q^N \neq 0$ for any integer $N \geq 1$.

Let $g = q(1 - e) - e \in \mathrm{End}(X \otimes X)$ and set

$$g_i = 1_{X^{i-1}} \otimes g \otimes 1 = q - (1 + q)e_i, \quad i = 1, 2, \ldots.$$

From the above relations for $e_1$ and $e_2$, the family $\{g_i\}_{i \geq i}$ is the generator of the Hecke algebra:

$$\begin{aligned} g_i g_{i+1} g_i = g_{i+1} g_i g_{i+1}, &\qquad i \geq 1, \\ g_i g_j = g_j g_i, &\qquad |i - j| \geq 2, \\ (g_i - q)(g_i + 1) = 0, &\qquad i \geq 1. \end{aligned}$$

Let $T = g_2 g_1 \in \mathrm{End}(X \otimes X \otimes X)$. Then we have

$$T(1_X \otimes \delta) = r(\delta \otimes 1_X) \qquad \text{for } \delta \in \mathrm{Hom}(I, X \otimes X),$$

where $r \in \mathbb{C}$ satisfies $r^2 = q^3$ (it is claimed in [12, Proposition 5.1] that $r^2 = 1$ with $N = 2$ but this should be corrected to $r^N = q^{N(N+1)/2}$). Note that, letting $t = r/q$, we have $q = t^2$, $r = t^3$ and, if $q$ is replaced with $q^{-1}$, then $r$ and therefore $t$ are changed into their inverses.

The tensor category is now completely classified by the complex number $t$ (or the pair $(q, r)$): Given a complex number $t$ satisfying $1 + t^2 + t^4 + \cdots + t^{2n} \neq 0$ for $n = 1, 2, \ldots$, there exists a unique rigid semisimple tensor category $\mathcal{C}(t)$ with the fusion rule governed by $S$. Two tensor categories $\mathcal{C}(t)$ and $\mathcal{C}(t')$ are isomorphic as monoidal categories if and only if either $t = t'$ or $t^{-1} = t'$.

These are in fact the Tannaka dual of the quantum group $SL_t(2, \mathbb{C})$: The tensor category $\mathcal{C}(t)$ is isomorphic to the one generated by the fundamental representation $V$ of the universal enveloping algebra $U_t$.

Recall the Hopf algebra $U_t$ is defined by the following relations on generators $\{K, E, F\}$ (see [10] for example):

$$KEK^{-1} = t^2 E, \quad KFK^{-1} = t^{-2} F, \quad [E, F] = \frac{K - K^{-1}}{t - t^{-1}},$$
$$\Delta(K) = K \otimes K, \quad \Delta(E) = E \otimes K + 1 \otimes E, \quad \Delta(F) = F \otimes 1 + K^{-1} \otimes F.$$

Note here that the relation

$$\tau(K) = K^{-1}, \quad \tau(E) = -t^{-1} E K^{-1}, \quad \tau(F) = -tKF$$

defines an antimultiplicative and anticomultiplicative involution of $U_t$.

If we define the fundamental representation of $U_t$ by $V = \mathbb{C}v_1 + \mathbb{C}v_2$ with

$$K = \begin{pmatrix} t & 0 \\ 0 & t^{-1} \end{pmatrix}, \quad E = \begin{pmatrix} 0 & 1 \\ 0 & 0 \end{pmatrix}, \quad F = \begin{pmatrix} 0 & 0 \\ 1 & 0 \end{pmatrix},$$

then we can easily check that the tensor category $\mathcal{R}$ generated by $V$ has the fusion rule $S$ with the invariant $(q, r) = (t^2, t^3)$, whence it is isomorphic to $\mathcal{C}(t)$.

With this model in hand, it is immediate to see that the tensor category $\mathcal{C}(t)$ admits a Frobenius duality ([20]): We first extend the Hopf algebra $U_t$ to $\widetilde{U}_t$ by adding a formal square root $K^{1/2}$ of $K$. By choosing a square root $t^{1/2}$ of $t$, we can extend $U_t$-modules in $\mathcal{R}$ to $\widetilde{U}_t$-modules without modifying the monoidal structure.

Given a $\widetilde{U}_t$-module $V$, the dual vector space $V^*$ is again a $\widetilde{U}_t$-module by

$$\langle xv^*, v \rangle = \langle v^*, \tau(x) v \rangle$$

and, given a $\widetilde{U}_t$-linear map $V \to W$, the usual transposed map ${}^t f$ is also $\widetilde{U}_t$-linear. Moreover, the obvious identification of the second dual $V^{**}$ with $V$ is $\widetilde{U}_t$-linear as $\tau$ being involutive.

Now, together with the involution explained so far, it is easy to see that the $\widetilde{U}_t$-linear map $\epsilon_V : V \otimes V^* \to \mathbb{C}$ defined by

$$\epsilon_V(v \otimes v^*) = \langle K^{1/2} v, v^* \rangle$$

for a $\widetilde{U}_t$-module $V$ gives a Frobenius duality in $\mathcal{R} \cong \mathcal{C}(t)$.

The quantum dimension of the fundamental representation $V$ is computed by

$$\epsilon_V {}^t\epsilon_V = \mathrm{tr}_V(K) = t + t^{-1}.$$

*Remark* . Given an integer $l \geq 3$, a similar analysis works without much difficulties for tensor categories of the truncated fusion rule of level $l$.

Given a finite family $\{\mathcal{C}(t_j)\}_{1\le j\le m}$ of tensor categories of type $A$ with their spectrum sets $\{S_j\}_{1\le j\le m}$, consider its free product $\mathcal{C} = \mathcal{C}(t_1)*\cdots*\mathcal{C}(t_m)$. The spectrum $S$ of $\mathcal{C}$ is then the free product $S = S_1 * \cdots * S_m$.

Letting $X_j$ be the fundamental generator of $\mathcal{C}(t_j)$ with the associated class $x_j \in S_j$, set $X = X_1 \otimes \cdots \otimes X_m$ in $\mathcal{C}$. Then we have the ascending sequence of semisimple algebras

$$\mathrm{End}(X) \subset \mathrm{End}(X \otimes X^*) \subset \mathrm{End}(X \otimes X^* \otimes X) \subset \dots .$$

We shall naturally identify the inductive limit algebra with the Bisch-Jones' planar algebra

$$FC(t_1 + t_1^{-1}, \dots, t_m + t_m^{-1}) = \bigcup_{n\ge 1} FC_n(t_1 + t_1^{-1}, \dots, t_m + t_m^{-1}),$$

where the coloring is specified by taking $\{x_1, \dots, x_m\}$ as the color set.

Given a periodic coloring

$$x_1x_2 \dots x_m x_m^* \dots x_1^* x_1 \dots x_m \dots$$

(although $x_j = x_j^*$, we dare to indicate $*$ depending on the parity of times of occurrence), we denote by $w_n$ its subword of the first $n$-colors for $n \ge 1$; $w_1 = x_1$, $w_2 = x_1x_2$, $w_{m+1} = x_1 \dots x_m x_m^*$ and so on.

According to [1], we introduce the algebra $A_n$ of planar string diagrams with both of the top and the bottom vertices colored by $w_n$ and the evaluation parameters given by

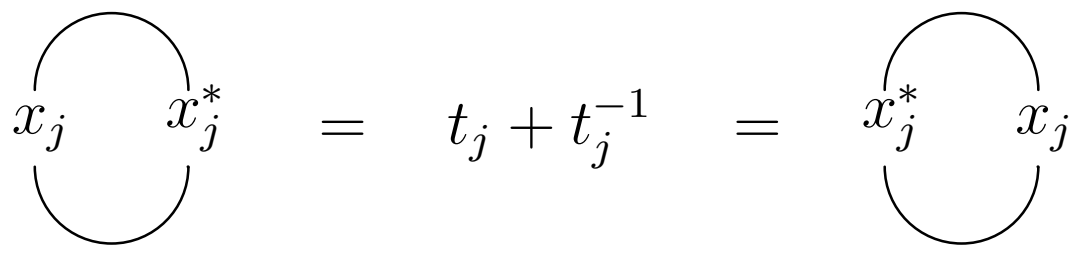

By definition, we have $A_{mn} = FC_n(t_1 + t_1^{-1}, \dots, t_m + t_m^{-1})$. In this sense, the algebras $A_n$ interpolate the Bisch-Jones planar algebras $FC_n$. As seen in [1, § 3.2], the dimension of $A_n$ (which is equal to the number of planar diagrams in $A_n$) is given by

$$\frac{l+1}{k(m+1)+l+1}\binom{k(m+1)+l+1}{k}$$

for $n = km + l$ with $0 \le l < m$.

By counting the parity of the number of end-points of strings, we have

**Lemma 4.1.** *In the planar algebra $A_n$, planar pairings among the source vertices (or the target vertices) occur only by coupling $x_j$ with $x_j^*$ (pairings of $x_j$ with $x_j$ or $x_j^*$ with $x_j^*$ being prohibited) while through strings connect $x_j$ with $x_j$ or $x_j^*$ with $x_j^*$.*

Let $W_n$ be the associated object of the word $w_n$ in $\mathcal{C} = \mathcal{C}_1 * \cdots * \mathcal{C}_m$: $W_1 = X_1$, $W_2 = X_1 \otimes X_2$, $W_{m+1} = X_1 \otimes X_m \otimes X_m^*$ and son on.

By the obvious imbedding $\mathrm{End}(W_n) = \mathrm{End}(W_n) \otimes 1 \subset \mathrm{End}(W_{n+1})$, we have the ascending sequence of semisimple algebras

$$\mathrm{End}(W_m) \subset \mathrm{End}(W_{m+1}) \subset \mathrm{End}(W_{m+2}) \subset \dots .$$

By the above lemma and hook identities for Frobenius duality, we can introduce an algebra homomorphism

$$\phi : \bigcup_{n\geq 1} A_n \to \bigcup_{n \geq 1} \mathrm{End}(W_n)$$

so that

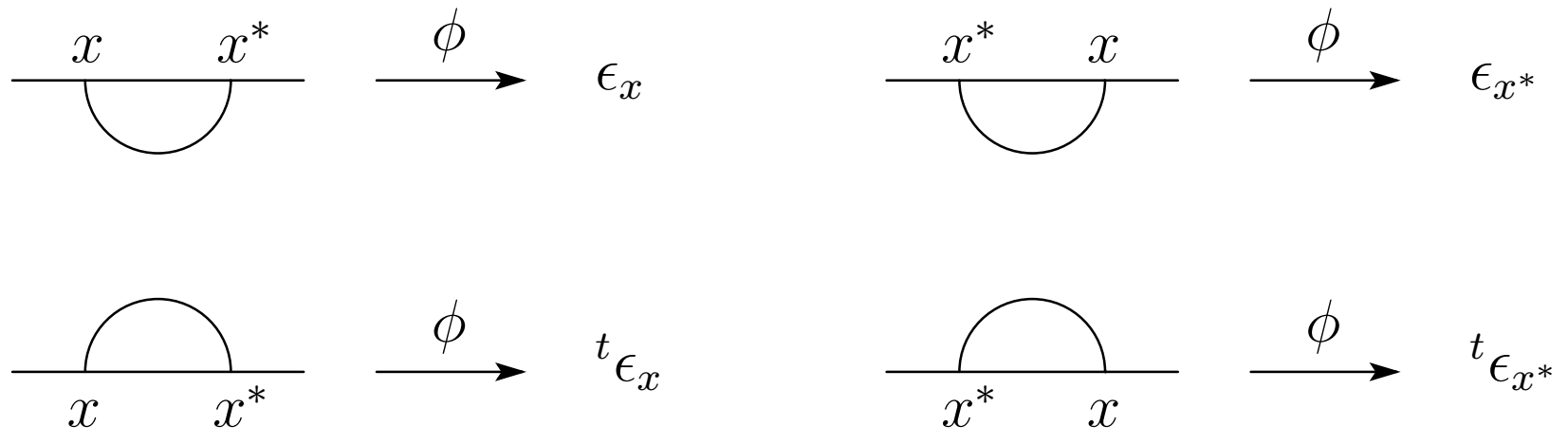

(diagrams stream up to down).

Given a simple object $s \in S = S_1 * \cdots * S_m$, the vector space $\mathrm{Hom}(s, W_n)$ with the obvious action of the algebra $\mathrm{End}(W_n)$ is a simple $\mathrm{End}(W_n)$-module whenever it is non-trivial.

In what follows, we shall use the letter $y$ to indicate one of $x_j$ $(1 \leq j \leq m)$ and define simple objects $\{y_k\}_{k\geq 0}$ inductively so that $y_0 = 1$, $y_1 = y$ and $y_k y = y y_k = y_{k-1} + y_{k+1}$.

We now introduce a one-to-one correspondance $s \leftrightarrow \sigma$ between elements in $S$ and words of the letter $\{x_1, \dots, x_m\}$, which is defined inductively so that (i) $s = \sigma$ if $|s| = 1$ and (ii) $\sigma = \sigma' y^k$ if $s = s' y_k$ with $s \| y_k$ and $s' \leftrightarrow \sigma'$, where $|s|$ denotes the length of $s \in S$.

Given a word $\sigma$ and an integer $n \geq 1$, let $V_\sigma^{(n)}$ be the $A_n$-module of the middle pattern $\sigma$ defined by [1, Definition 3.1.16]. Recall that $V_\sigma^{(n)}$ is the free vector space over the set of colored planar diagrams from $\sigma$ to $w_n$ with no self-couplings among vertices in $\sigma$ and the action of $A_n$ is defined by the obvious composition of diagrams with the convention that, if composed diagrams get out of $V_\sigma^{(n)}$, then it is set to be zero.

We now introduce a linear map $\Phi : V_\sigma^{(n)} \to \mathrm{Hom}(s, W_n)$ in the following way: By replacing self-coupling inside $w_n$ by ${}^t\epsilon$'s, each planar diagram in $V_\sigma^{(n)}$ associates a morphism in $\mathrm{Hom}(\otimes\sigma, W_n)$, where $\otimes\sigma$ denotes the tensor product of objects appearing in $\sigma$ with $*$ placed according to the parity of times of occurrence (for example, $\otimes\sigma =$ $X_1 \otimes X_1^* \otimes X_1 \otimes X_3 \otimes X_3^*$ for $\sigma = x_1x_1x_1x_3x_3$). Now, for each $k \geq 1$ and each $y$ in $\{x_1, \dots, x_m\}$, choose a non-trivial morphism

$$y_k \to y^{\otimes k} = y \otimes y^* \otimes y \otimes \dots$$

once for all (which is determined up to scalar multiplication) and form the tensor product for factors appearing in $\otimes\sigma$, resulting a morphism $s \to \otimes\sigma$. By taking the composition with this, we obtain the linear map

$$\Phi : V_\sigma^{(n)} \longrightarrow \mathrm{Hom}(\otimes\sigma, W_n) \longrightarrow \mathrm{Hom}(s, W_n).$$

**Lemma 4.2.** *For $a \in A_n$, we have the commutative diagram*

$$\begin{array}{ccc} V_\sigma^{(n)} & \longrightarrow & \mathrm{Hom}(s, W_n) \\ {\scriptstyle\sigma}\big\downarrow & & \big\downarrow{\scriptstyle\phi(a)} \\ V_\sigma^{(n)} & \longrightarrow & \mathrm{Hom}(s, W_n) \end{array} .$$

*Proof.* If $a$ and $v$ are planar diagrams in $A_n$ and $V_\sigma^{(n)}$ respectively with the composed diagram $av$ out of $V_\sigma^{(n)}$, then the image of $av$ in $\mathrm{Hom}(\otimes\sigma, W_n)$ is factored into the form $f \circ g$, where $g \in \mathrm{Hom}(\otimes\sigma, \otimes\sigma')$ and $f \in \mathrm{Hom}(\otimes\sigma', W_n)$ with $\sigma'$ a strict subsequence of $\sigma$. Then, by an easy induction, we see that $\mathrm{Hom}(s, \otimes\sigma') = 0$ (the highest spin part $y_k$ in $y^{\otimes k}$ does not appear in $y^{\otimes j}$ for $j < k$), whence $av$ is in the kernel of the linear map $V_\sigma^{(n)} \to \mathrm{Hom}(s, W_n)$. □

**Lemma 4.3.** *The linear map $\Phi : V_\sigma^{(n)} \to \mathrm{Hom}(s, W_n)$ is an isomorphism for any $\sigma$ and any $n \geq 1$.*

*The representation of $A_n$ on $\mathrm{Hom}(s, W_n)$ via $\phi$ is irreducible and inequivalent for different $s$ (as long as $\mathrm{Hom}(s, W_n) \neq 0$).*

*Proof.* We prove the assertion by an induction on $n$. Assume that the statement is true for $n$. Consider the $\mathrm{End}(W_{n+1})$-module $\mathrm{Hom}(s_{n+1}, W_{n+1})$ with $s_{n+1} \in S$ and let us begin with looking at how it is decomposed into simple components when restricted to $\mathrm{End}(W_n)$. Let $w_{n+1} = w_n y$ with $y \in \{x_1, \dots, x_m\}$ and write $W_{n+1} = W_n \otimes y$. By the obvious isomorphism

$$\mathrm{Hom}(s_{n+1}, W_{n+1}) \cong \bigoplus_{s \in S} \mathrm{Hom}(s, W_n) \otimes \mathrm{Hom}(s_{n+1}, s \otimes y),$$

together with the Frobenius reciprocity

$$\mathrm{Hom}(s_{n+1}, s \otimes y) \cong \mathrm{Hom}(s_{n+1} \otimes y^*, s),$$

we need to decompose $s_{n+1}y$: If $s_{n+1}\|y$, $s_{n+1}y \in S$ and $s$ is forced to be $s_{n+1}y$. Then we have the isomorphism

$$\mathrm{Hom}(s_{n+1}y, W_n) \ni w \mapsto (w \otimes 1_y)(1 \otimes {}^t\epsilon_{y^*}) \in \mathrm{Hom}(s_{n+1}, W_{n+1})$$

and the $A_n$-module $\mathrm{Hom}(s_{n+1}, W_{n+1})$ is equivalent to the $A_n$-module $\mathrm{Hom}(s_{n+1}y, W_n)$, which is irreducible by the induction hypothesis. Thus the $A_{n+1}$-module $\mathrm{Hom}(s_{n+1}, W_{n+1})$ is irreducible as well.

Next we consider the case that $s_{n+1}|y$ with $s_{n+1} = sy_k$ ($k \geq 1$). From $s_{n+1}y = sy_{k-1}+sy_{k+1}$, we have two simple components $\mathrm{Hom}(sy_{k-1}, W_n)$ and $\mathrm{Hom}(sy_{k+1}, W_n)$, which are mapped into $\mathrm{Hom}(s_{n+1}, W_{n+1})$ by

$$\mathrm{Hom}(s_n, W_n) \ni w \mapsto (w \otimes 1_y)(1_s \otimes f) \in \mathrm{Hom}(s_{n+1}, W_{n+1})$$

with $s_n = sy_{k-1}$ or $s_n = sy_{k+1}$ and $0 \neq f$ in $\mathrm{Hom}(y_k, y_{k-1} \otimes y)$ or $\mathrm{Hom}(y_k, y_{k+1} \otimes y)$ respectively.

The $A_n$-module $\mathrm{Hom}(s_{n+1}, W_{n+1})$ is then decomposed into a direct sum of two inequivalent $A_n$-modules $\mathrm{Hom}(sy_{k-1}, W_n)$ and $\mathrm{Hom}(sy_{k+1}, W_n)$ by the induction hypothesis.

To see the irreducibility of the $A_{n+1}$-module $\mathrm{Hom}(s_{n+1}, W_{n+1})$, we choose a planar diagram $v$ in $V_{\sigma y^{k+1}}$ and let $w$ be its image in $\mathrm{Hom}(sy_{k+1}, W_n)$, which corresponds to

$$\widetilde{w} = (w \otimes 1_y)(1_s \otimes f) \quad \text{with } f : y_k \to y_{k+1} \otimes y$$

in $\mathrm{Hom}(s_{n+1}, W_{n+1})$. Let $\widetilde{v} \in V^{(n+1)}_{\sigma y^k}$ be the diagram obtained by taking the contraction of $v \otimes 1_y \in V^{(n+1)}_{\sigma y^{k+2}}$ for two through strings colored by $y$ (and $y^*$) at the right end and set $a = \widetilde{v}\widetilde{v}^* \in A_{n+1}$ (Fig. 5). Then we have $\phi(a)\widetilde{w} = (\prod \mathrm{O})\, \Phi(\widetilde{v}) \circ h$ with $h = (1_{y_k} \otimes \epsilon_{y^*})(g \otimes 1_y)f$ and $g : y_{k+1} \to y_k \otimes y^*$ (Fig. 8). By Frobenius reciprocity, $h \in \mathrm{End}(y_k)$ has the expression $(1_{y_k} \otimes \epsilon_{y^*})(g\widetilde{f} \otimes 1_y)(1_{y_k} \otimes {}^t\epsilon_{y^*})$ and its quantum trace is calculated by

$$\epsilon_{y_k}(h \otimes 1)^t\epsilon_{y_k} = \epsilon_{y_k y^*}(g\widetilde{f} \otimes 1)^t\epsilon_{y_k y^*} = \epsilon_{y_{k+1}}(\widetilde{f}g \otimes 1)^t\epsilon_{y_{k+1}} = c\frac{t_j^{k+2} - t_j^{-k-2}}{t_j - t_j^{-1}}$$

if $\widetilde{f}g = c1_{y_{k+1}}$ and $y = x_j$, which does not vanish. Thus $\phi(a)\widetilde{w}$ is a non-zero constant multiple of $\Phi(\widetilde{v})$, which belongs to the image of $\mathrm{Hom}(sy_{k-1}, W_n)$ in $\mathrm{Hom}(s_{n+1}, W_{n+1})$; $\mathrm{Hom}(sy_{k+1}, W_n)$ and $\mathrm{Hom}(sy_{k-1}, W_n)$ are mixed up by the action $\phi(a)$ on $\mathrm{Hom}(s_{n+1}, W_{n+1})$. This, together with the fact $\phi(A_n)' \cong \mathbb{C} \oplus \mathbb{C}$, shows that $\phi(A_{n+1})' = \mathbb{C}$, proving the irreducibility of the $A_{n+1}$-module $\mathrm{Hom}(s_{n+1}, W_{n+1})$.

Since the image of $V^{(n+1)}_{\sigma y^k}$ in $\mathrm{Hom}(sy_k, W_{n+1})$ is clearly non-trivial, the irreducibility particularly shows the surjectivity of $\Phi$.

On the other hand, we know the exact sequence of $A_n$-modules

$$0 \longrightarrow V^{(n)}_{\sigma y^{k-1}} \longrightarrow V^{(n+1)}_{\sigma y^k} \longrightarrow V^{(n)}_{\sigma y^{k+1}} \longrightarrow 0$$

(see the proof of [1, Theorem 3.2.1]) and the induction hypothesis implies

$$\begin{aligned}\dim V^{(n+1)}_{\sigma y^k} &= \dim V^{(n)}_{\sigma y^{k-1}} + \dim V^{(n)}_{\sigma y^{k+1}} \\ &= \dim \mathrm{Hom}(sy_{k-1}, W_n) + \dim \mathrm{Hom}(sy_{k+1}, W_n) \\ &= \dim \mathrm{Hom}(sy_k, W_{n+1}),\end{aligned}$$

proving the bijectivity of the map $\Phi : V^{(n+1)}_{\sigma y^k} \to \mathrm{Hom}(sy_k, W_{n+1})$.

The inequivalence of $A_{n+1}$-modules $\mathrm{Hom}(s_{n+1}, W_{n+1})$ for different $s_{n+1}$ is now obvious because $s_{n+1} \in S$ can be recovered by checking the irreducible decomposition of the restriction to $A_n$. ☐

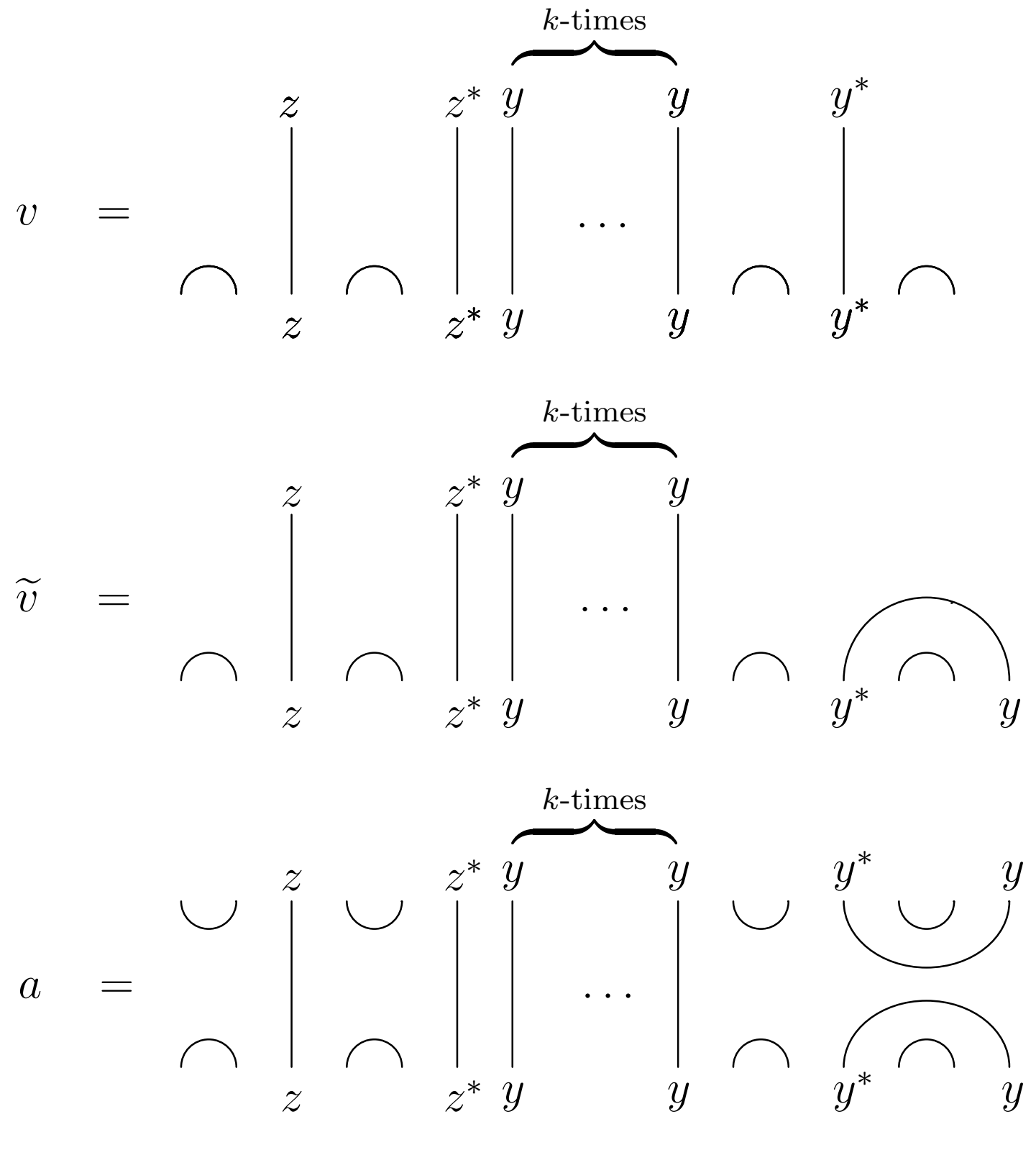

Figure 5.

The following should be compared with [1, Theorem 6.1.2].

**Proposition 4.4.** *The map $\phi$ is bijective: The Bisch-Jones' planar algebra is identified with the inductive limit of the ascending sequence*

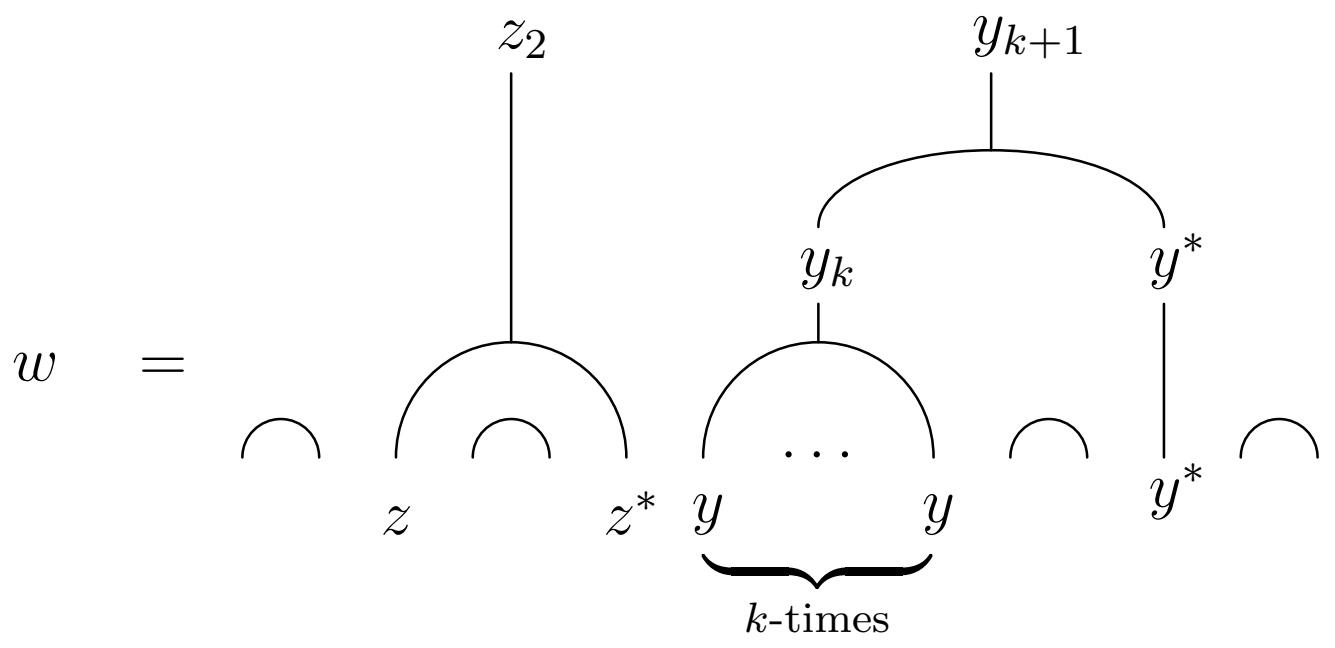

FIGURE 6.

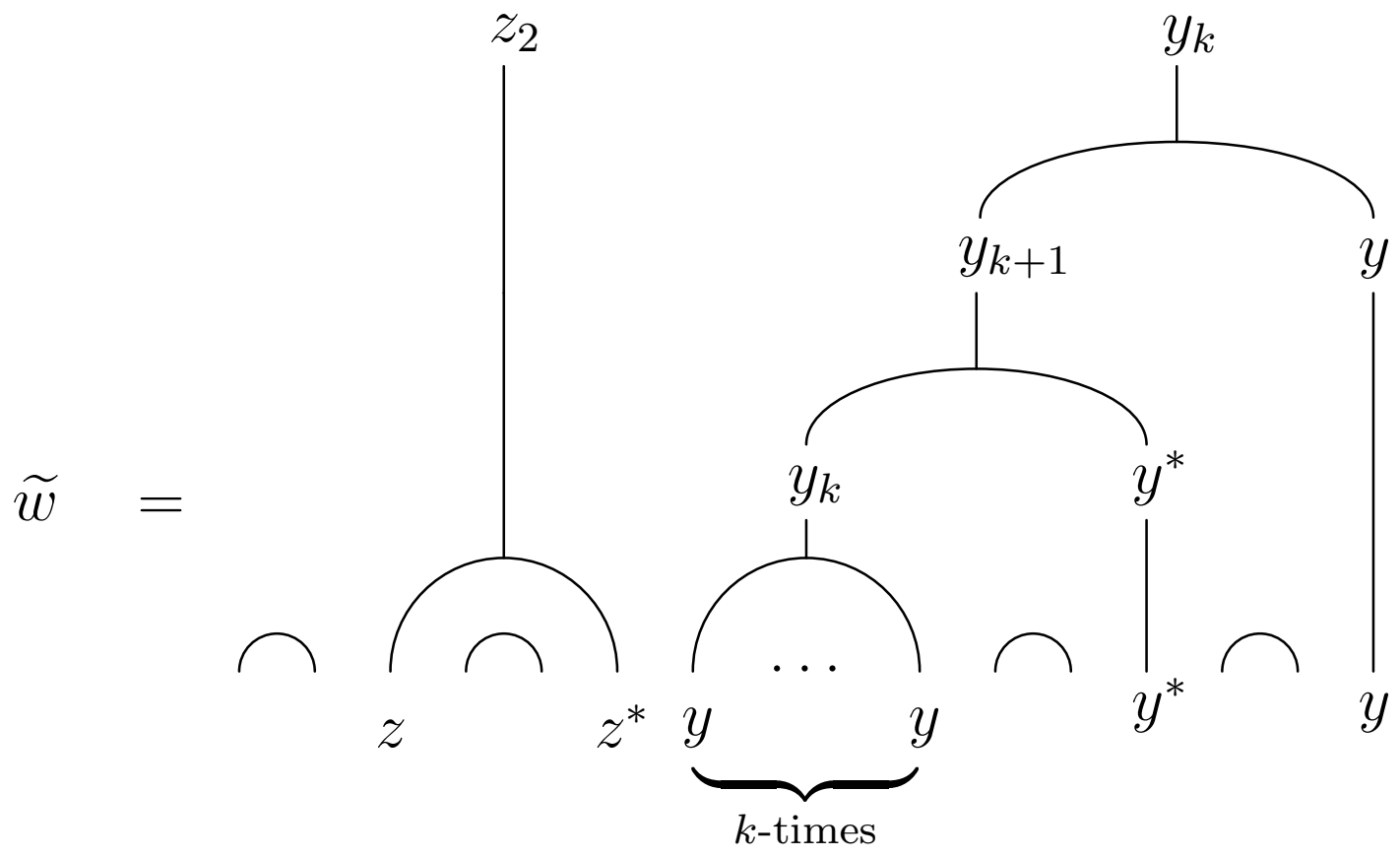

FIGURE 7.

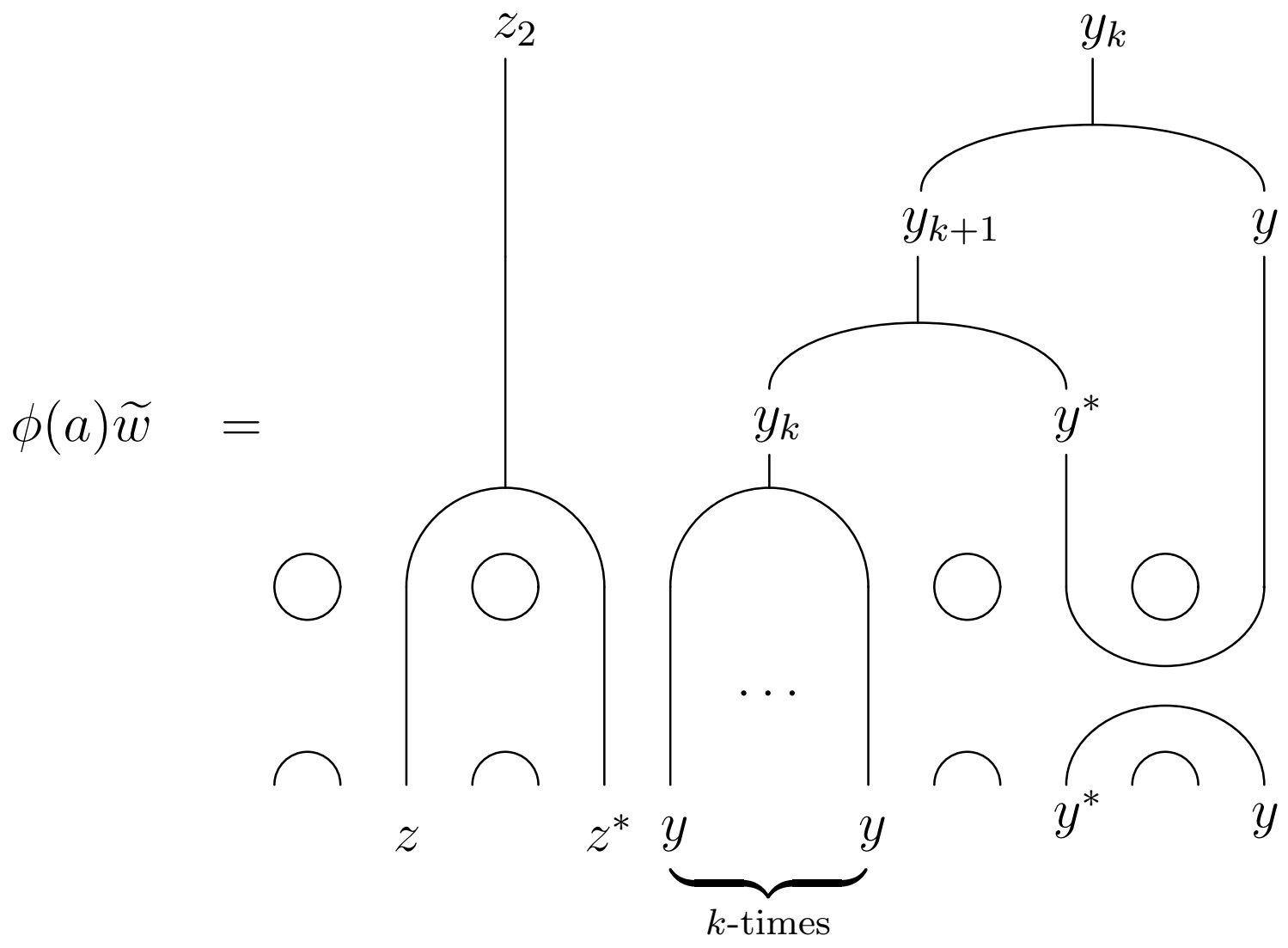

FIGURE 8.

*of algebras associated to the simple object $x = x_1 x_2 \dots x_n$ in the free product tensor category $\mathcal{C}$.*

*Proof.* By the previous lemma, we see that $\phi : A_n \to \mathrm{End}(W_n)$ is surjective because $\mathrm{End}(W_n)$ is semisimple and inequivalent irreducible representations of $\mathrm{End}(W_n)$ gives rise to inequivalent irreducible representations of $A_n$. (Use the double commutant theorem for the semisimple $A_n$-module $\bigoplus_{s \in S} \mathrm{Hom}(s, W_n)$.)

From the analysis in [1], we know that

$$\begin{aligned} \dim A_n &= \sharp\{\text{planar diagrams in } A_n\} \\ &= \text{the dimension of } \mathrm{End}(W_n) \text{ calculated by the fusion rule} \\ &= \dim \mathrm{End}(W_n), \end{aligned}$$

which is checked for generic evaluation parameters but the formula itself clearly holds without restrictions.

Since $\phi_n$ is surjective, the equality of dimensions shows that it is in fact bijective. □

**Corollary 4.5.** *The Bisch-Jones' planar algebra $FC_n(a_1, \dots, a_m)$ is semisimple for all $n \geq 1$ if and only if none of $a_1, \dots, a_m$ belongs to the set*

$$\{2\cos(\pi r); r \in \mathbb{Q} \setminus \mathbb{Z}\}.$$

*Proof.* With the choice $a = t + t^{-1}$, the condition $1 + t^2 + t^4 + \cdots + t^{2n} \neq 0$ for any $n \geq 1$ is equivalent to require that $a$ does not belong to the set specified above.

Conversely assume that the Fuss-Catalan algebra $FC_n$ is semisimple for any $n \geq 1$ and look at the coloring $x_j$ $(1 \leq j \leq m)$. In the planar algebra $A_{mn} = FC_n$ with $n \geq 1$ an odd integer, taking all the possible pairings of the form $\epsilon_{x_{j+1} \dots x_m}$ or $\epsilon_{x_{j-1} \dots x_1}$ (multiplied by $a_{j+1}^{-1} \dots a_m^{-1}$ or $a_1^{-1} \dots a_{j-1}^{-1}$ respectively), we obtain an idempotent $e \in A_{mn}$ of middle pattern $x_1 \dots x_{j-1} x_j^n x_{j+1} \dots x_m$.

By the choice, the reduced algebra $eA_{mn}e$ is isomorphic to the Temperley-Lieb algebra, which is semisimple by our assumption. It is well-known ([18, 6]) that the Temperley-Lieb algebra of $n$ strings ([11]) is semisimple for all odd $n \geq 1$ (if and) only if the parameter $a_j$ is out of the range of numbers in question. □

The following should be compared with [1, Theorem 3.3.4] (the explicit formula for minimal projections is not needed here).

**Corollary 4.6.** *If the quantum (non-normalized) trace is defined on $FC_n$ by multiplying evaluation parameters of loops obtained after closing diagrams, then the value $d_\sigma$ of the trace on a minimal projection* *corresponding to the irreducible representation* $V_\sigma^{(mn)}$ *($\sigma$ being a word of letters $x_1, \dots, x_m$) is inductively calculated by the rule $d_\sigma = t_j + t_j^{-1}$ for $\sigma = x_j$ and*

$$d_\sigma = d_{\sigma'} \frac{t_j^{k+1} - t_j^{-k-1}}{t_j - t_j^{-1}}$$

*if $\sigma = \sigma' x_j^k$ with $x_j$ different from the last letter of $\sigma'$.*

*Moreover, we have the identity*

$$(t_1 + t_1^{-1})^n \dots (t_m + t_m^{-1})^n = \sum_\sigma d_\sigma \dim V_\sigma^{(mn)}.$$

*Note that* $\dim V_\sigma^{(mn)}$ *is equal to the number of possible planar diagrams inside.*

## Appendix A.

Here we shall describe the coherence theorem on triangulated vector spaces of a labeled polygon.

Assume that we are given a fusion rule set $S$ (no need to have the unit nor the involution here), triangular vector spaces $\left\{ \begin{bmatrix} x\,y \\ z \end{bmatrix} \right\}$ and associativity transformations

$$\bigoplus_{x_{12} \in S} \begin{bmatrix} x_1 \; x_2 \\ x_{12} \end{bmatrix} \otimes \begin{bmatrix} x_{12} \; x_3 \\ x_0 \end{bmatrix} \longrightarrow \bigoplus_{x_{23} \in S} \begin{bmatrix} x_2 \; x_3 \\ x_{23} \end{bmatrix} \otimes \begin{bmatrix} x_1 \; x_{23} \\ x_0 \end{bmatrix}$$

Let $P$ be a (convex) polygon with a distinguised edge placed at the bottom. By a labeling of $P$, we shall mean an assignment of elements to edges in $P$. Given a triangulation $T$ of a labeled polygon $P$, we can associate the vector space $[T]$ by taking the summation of possible tensor products of triangular vector spaces.

If one applies associativity transformations locally, we obtain isomorphisms among vector spaces $[T]$ for variety of choices of $T$ and here comes out the problem of coherence: Let $A$, $B$ be two triangulations of $P$ and $A = S_1, S_2, \cdots, S_m = B$, $A = T_1, T_2, \cdots, T_n = B$ be sequences of triangulations such that $[S_{j-1}]$ and $[S_j]$ (resp. $[T_{j-1}]$ and $[T_j]$) can be related by a single associativity transformation. We then obtain two isomorphisms between $[A]$ and $[B]$ as succesive applications of associativity transformations according to the histories $(S_1, \cdots, S_m)$ and $(T_1, \cdots, T_n)$. The problem of coherence is then whether we can deduce the equality of these isomorphisms or not. The problem is obviously reduced to the case of loops, i.e., $A = B$, and it suffices to fix the starting configuration $T$ once for all.

**Proposition A.1** (Coherence Theorem)**.** *If associativity transformations satisfy the coherence for pentagons, then they are coherent for general polygons.*

We need some terminologies. Let $P$ be a polygon and denote by $\mathcal{T}$ the set of triangulations of $P$. The set $\mathcal{T}$ is made into the vertex set of a graph by joining two triangulations with a single edge if they are related by an associativity transformation. We choose a distinguished vertex $\bullet$ in $P$ and define the triangulation $T_0$ by drawing all the possible diagonal lines passing through the distinguished vertex.

The length $l(T)$ of a triangulation $T$ in $\mathcal{T}$ is, by definition, the number of diagonal lines in $T$ penetrating the distinguished vertex. $T_0$ is the unique triangulation of maximal length. A directed path in $\mathcal{T}$ defines the accompanied isomorphism of vector spaces by succesive applications of (amplified) associativity transformations. A directed path in $\mathcal{T}$ starting at $T_0$ is called a **short-cut** if the length is strictly decreasing along the path. It is easy to see that any triangulation has a short-cut.

**Lemma A.2.** *Let $T$ be a triangulation in $\mathcal{T}$. Then all the short-cuts from $T_0$ to $T$ give the same isomorphism.*

*Proof.* We prove by a (reverse) induction on the length of $T$. If $l(T) = 0$, $T$ is adjacent to the unique triangulation of length 1 and the problem is reduced to that for $l(T) = 1$.

For $l(T) \geq 1$, $T$ contains a diagonal line $L$ passing through the vertex $\bullet$. Let $P'$ and $P''$ be subpolygons of $P$ separated by $L$ with the induced triangulations $T'$ and $T''$ respectively ($T = T' \times T''$). We may assume that $P''$ contains the bottom edge $B$ without loss of generality (see Fig. 9). Since any short-cut to $T$ does not change the line $L$, associativity transformations in a short-cut $T_0, T_1, \dots, T_n, T$ $(n \geq 0)$ are operations on $P'$ or $P''$. Thus, gathering these into two groups, we obtain short-cuts $T'_1, \dots, T'_{n'}$ and $T''_1, \dots, T''_{n''}$ for $T'$ and $T''$ respectively. By associativity of tensor products of linear maps, the isomorphism $[T_0] \to [T]$ specified by the short-cut $\{T_j\}$ is equal to the one given by the short-cut $\{T'_1 \times T''_0, \dots, T'_{n'} \times T''_0, T' \times T''_0, T' \times T''_1, \dots, T' \times T''_{n''}\}$.

In this way, the problem is reduced to those for $T'$ and $T''$, which obviously have smaller lengths than $T$. □

We now prove the coherence theorem in the form that any closed path

$$\gamma : T_0 \to T_1 \to T_2 \to \cdots \to T_{k-1} \to T_k \to T_0$$

gives the identity transformation. We first connect each triangulation $T_i$ by a short-cut (from $T_0$). Then it suffices to show the commutativity

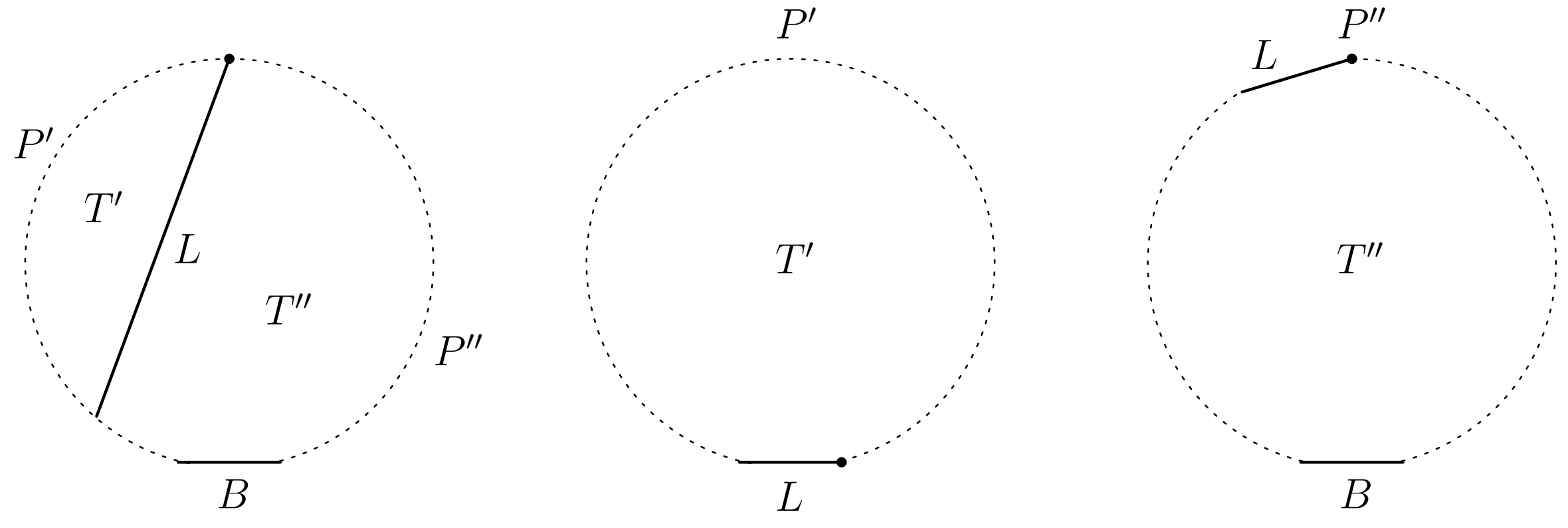

FIGURE 9.

of the following triangular diagram given by associativity transformations:

$$\begin{array}{ccc} T_j & \longrightarrow & T_{j+1} \\ {\scriptstyle \text{short-cut}}\uparrow & & \uparrow{\scriptstyle \text{short-cut}} \\ T_0 & = & T_0 \end{array}$$

If $|l(T_j) - l(T_{j+1})| = 1$, say $l(T_{j+1}) = l(T_j) - 1$, then comparing two short-cuts $T_0 \to T_j \to T_{j+1}$ and $T_0 \to T_{j+1}$, the previous lemma ensures the assertion.

Assume on the contrary that $l(T_j) = l(T_{j+1}) = l$. Recall here that we can use any short-cuts to prove the commutativity of the diagram by the previous lemma. The problem in question will then be worked out by an induction on the common length $l$. We consider two cases: If we can find a square which does not touch the operation $T_j \to T_{j+1}$ and contains the distinguished vertex $\bullet$ with $T_j$ and $T_{j+1}$ not including the diagonal lines of the square passing through $\bullet$ (Fig. 10), then we define the triangulations $S_j$ and $S_{j+1}$ in $\mathcal{T}$ of length $l+1$ by changing the diagonal line of the square in $T_j$ and $T_{j+1}$ (Fig. 11), which yield the commutative diagram

$$\begin{array}{ccc} T_j & \longrightarrow & T_{j+1} \\ \uparrow & & \uparrow \\ S_j & \longrightarrow & S_{j+1} \end{array}$$

and hence the problem is reduced to the lower level $S_j$ and $S_{j+1}$.

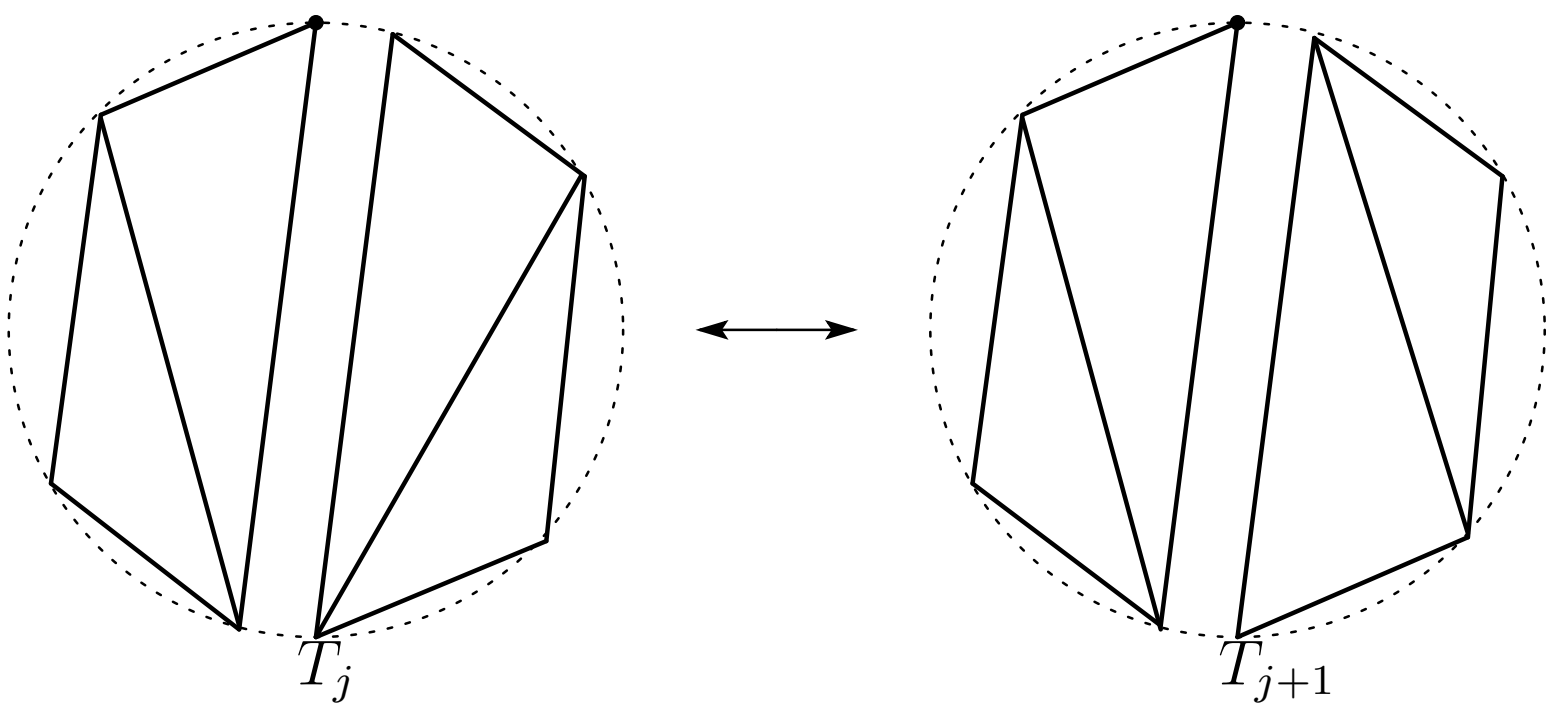

FIGURE 10.

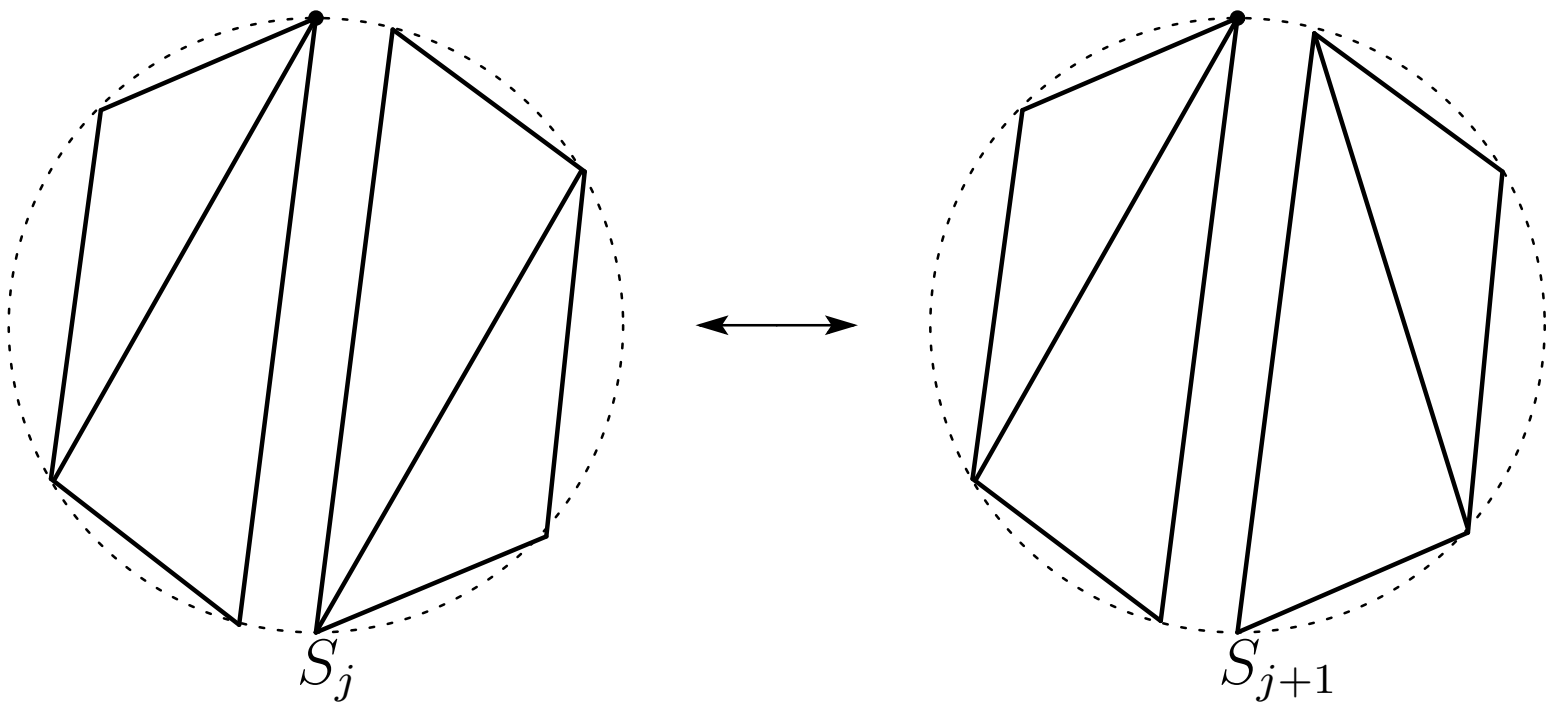

FIGURE 11.

The remaining is the case that, for any square containing $\bullet$ with the diagonal line not crossing $\bullet$, the operation $T_j \iff T_{j+1}$ must change one of edges in this square. Since $l(T_j) = l(T_{j+1})$, edges passing though $\bullet$ cannot be changed. Thus the possible situation is described as Fig. 12.

Now we define $S_j$, $S_{j+1}$, $S \in \mathcal{T}$ by Fig. 13 so that they, together with $T_j$ and $T_{j+1}$, form a pentagonal diagram

$$\begin{array}{ccccc} S & \longrightarrow & S_j & \longrightarrow & T_j \\ \| & & & & \downarrow \\ S & \longrightarrow & S_{j+1} & \longrightarrow & T_{j+1}. \end{array}$$

Thus, by taking a short-cut to $S$ and then applying the pentagonal relation of associativity transformations, we complete the proof of the commutativity of the short-cut diagram.

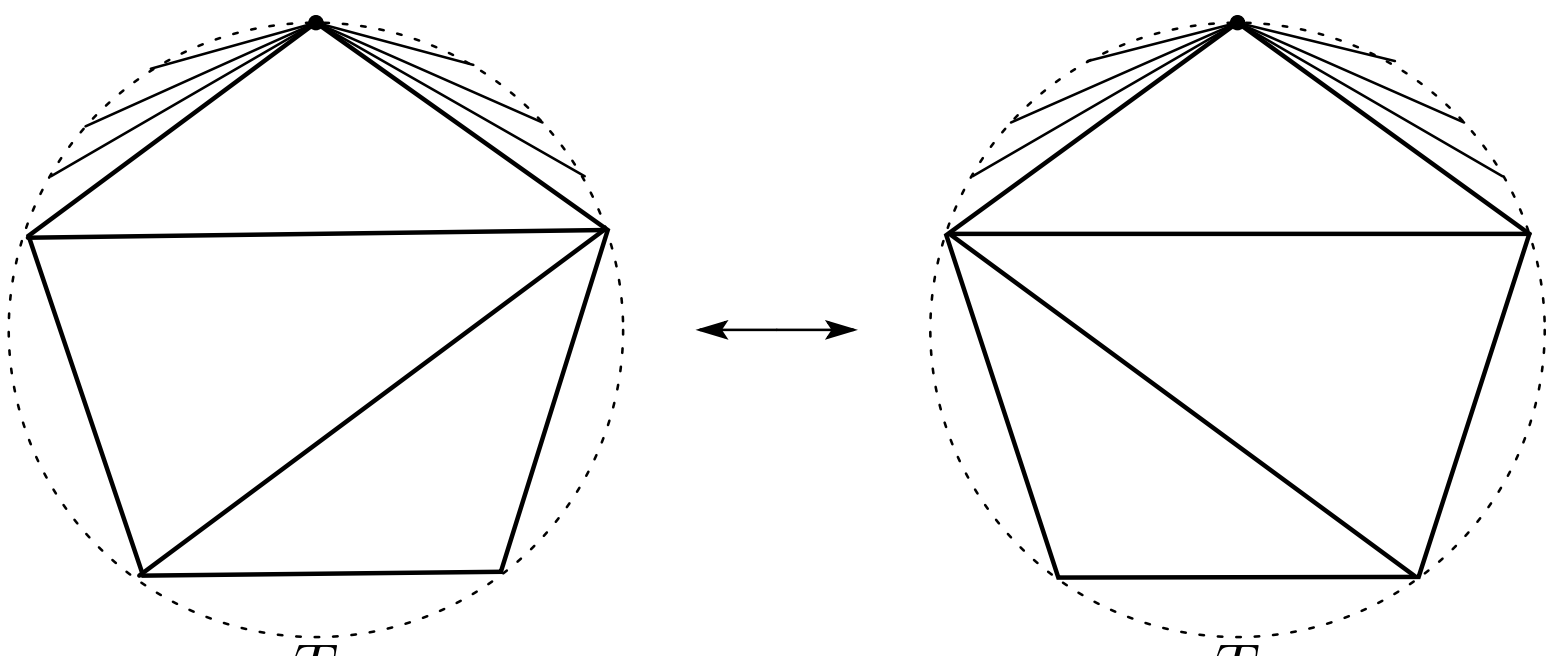

FIGURE 12.

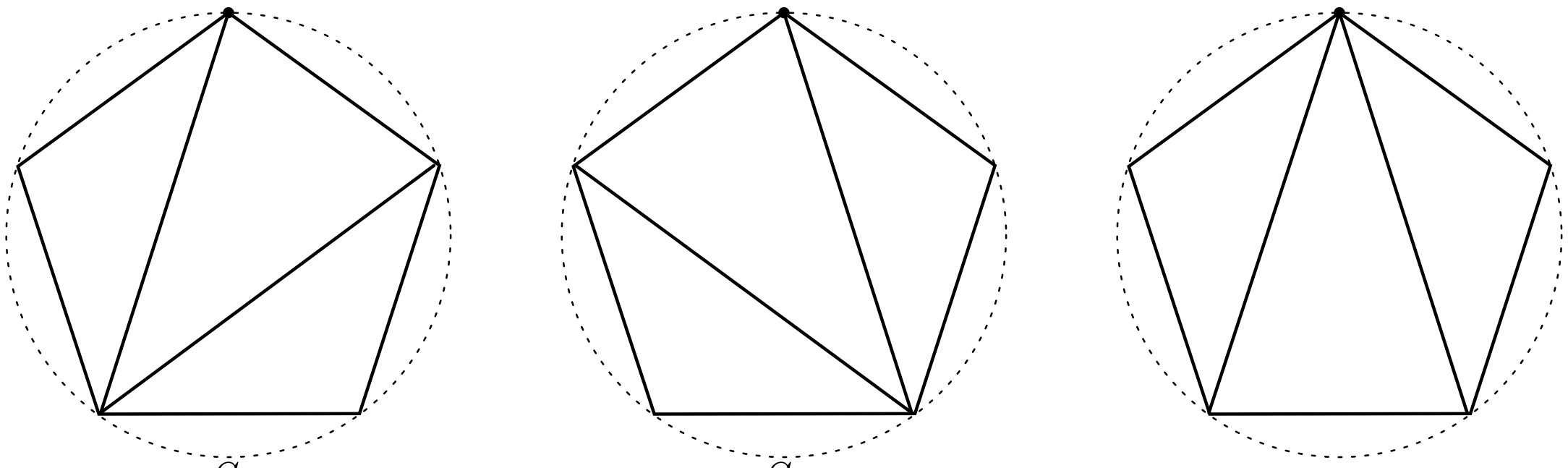

FIGURE 13.